\documentclass[a4paper,12pt]{article} %izvzela 'twoside' fleqn,

\usepackage{color}

\usepackage{amsfonts, amsmath, amsthm, amssymb}
\usepackage[T1]{fontenc}
\usepackage[cp1250]{inputenc}
\usepackage{xcolor}
\usepackage{graphicx}
\usepackage{amssymb}
\usepackage{amsmath}
\usepackage{mathptmx}
\usepackage{helvet}
\usepackage{courier}
\usepackage{txfonts}
\usepackage{tikz} 
\usetikzlibrary{arrows}
\usepackage{type1cm}

\usepackage{verbatim}

\usepackage{graphicx}
\usepackage{epsfig,amscd,amssymb,amsxtra,amsmath,amsthm}
\usepackage{type1cm}
\usepackage[T1]{fontenc}
\usepackage{graphics}
\usepackage[mathscr]{eucal}
\usepackage[all]{xy}
\usepackage{amsmath,amscd}

\newtheorem{theorem}{Theorem}[section]

\newtheorem{definition}[theorem]{Definition}
\newtheorem{lemma}[theorem]{Lemma}

\newtheorem{example}[theorem]{Example}

\newtheorem{corollary}[theorem]{Corollary}
\newtheorem{problem}[theorem]{Problem}
\newtheorem{observation}[theorem]{Observation}

\newcommand{\diam} {\mathop{\rm diam}\nolimits}

\newcommand{\Cl}  {\mathop{\rm Cl}\nolimits}

\begin{document}

\def\joinrel{\mkern-3mu}
\newcommand{\varproj}{\displaystyle \lim_{\multimapinv\joinrel-\joinrel-}}

\title{Quotients of dynamical systems and chaos on the Cantor fan}
\author{Iztok Bani\v c, Goran Erceg, Judy Kennedy and Van Nall}
\date{}

\maketitle

\begin{abstract}
Let $(X,f)$ be a dynamical system. Using an equivalence relation $\sim$ on $X$, we introduce the quotient $(X/_{\sim},f^{\star})$ of the dynamical system $(X,f)$. In the first part of the paper, we give new results about sensitive dependence on initial conditions of  $(X/_{\sim},f^{\star})$, transitivity of $(X/_{\sim},f^{\star})$, and periodic points in $(X/_{\sim},f^{\star})$. In the second part of the paper, we use these results to study chaotic functions on the Cantor fan. Explicitly, we study functions $f$ on the Cantor fan $C$ such that (1) $(C,f)$ is chaotic in the sense of Devaney, (2) $(C,f)$ is chaotic in the sense of Robinson but not  in the sense of Devaney, and (3) $(C,f)$  is chaotic in the sense of Knudsen but not in the sense of Devaney. We also study chaos on the Lelek fan.  
\end{abstract}
\-

%%%%%%%%%%%%%%%%%%%%%%%%%%%%%%%%%%%%%%%%%%%%%%%%%%%%%%%%%%%%%%%%%%%%%%%%%%%%%%%%%
%%% I N T R O D U C T I O N S
\section{Introduction}
A dynamical system is a pair $(X,f)$, where $X$ is a compact metric space and $f:X\rightarrow X$ is a continuous function. Recently, many interesting properties of dynamical systems have been obtained in such a way that from a given dynamical system $(X,f)$ a new dynamical system was created by defining an equivalence relation on $X$ and then, by transforming the function $f$ according to the defined equivalence relation;  see \cite{BE,BE1,BE2} for such examples. The following theorem is used as a main tool to obtain these results.
\begin{theorem}\label{kvocienti}
Let $X$ be a compact metric space, let $\sim$ be an equivalence relation on $X$, and  let $f:X\rightarrow X$ be a function such that for all $x,y\in X$,
$$
x\sim y  \Longleftrightarrow f(x)\sim f(y).
$$
 Then the following hold for $f^{\star}$, which is defined by $f^{\star}([x])=[f(x)]$ for any $x\in X$. 
\begin{enumerate}
\item\label{u1} $f^{\star}$ is a well-defined function from  $X/_{\sim}$ to $X/_{\sim}$. 
\item\label{u2} If $f$ is continuous, then $f^{\star}$ is continuous.
\item\label{u3} If $f$ is a homeomorphism, then $f^{\star}$ is a homeomorphism.
\item\label{u4} Suppose that $X/_{\sim}$ is metrizable. If $(X,f)$ is transitive, then $(X/_{\sim},f^{\star})$ is transitive.
\end{enumerate}
\end{theorem}
\begin{proof}
	See \cite[Theorem 3.4]{BE}. 
\end{proof}
We call the pair $(X/_{\sim},f^{\star})$ from Theorem \ref{kvocienti}, the quotient of the dynamical system $(X,f)$. Note that if $X/_{\sim}$ is metrizable and if  $(X,f)$ is a dynamical system, then also  $(X/_{\sim},f^{\star})$ is a dynamical system.
In the first part of the paper, we discuss the following problem, which is motivated by Theorem \ref{kvocienti}.
\begin{problem}
	Let $X$ be a compact metric space, let $\sim$ be an equivalence relation on $X$,  let $f:X\rightarrow X$ be a function such that for all $x,y\in X$,
$$
x\sim y  \Longleftrightarrow f(x)\sim f(y),
$$
and let $\mathcal P$ be a dynamical property. 
If $(X,f)$ is a dynamical system that has the property $\mathcal P$, does then the quotient $(X/_{\sim},f^{\star})$ also have the property $\mathcal P$ (if $X/_{\sim}$ is metrizable)?
\end{problem}
We study dynamical systems $(X,f)$ with property $\mathcal P$, where $\mathcal P$ is either sensitive dependence on initial conditions, transitivity, or the property that the set $\mathcal P(f)$ of periodic points in $(X,f)$ is dense in $X$.  In particular, we give an example of a dynamical system $(X,f)$ which has sensitive dependence on initial conditions but the quotient $(X/_{\sim},f^{\star})$ does not; see Example \ref{buda}. We show in Theorem \ref{kvocki} that under some minor additional conditions, $(X/_{\sim},f^{\star})$ always has sensitive dependence on initial conditions, if $(X,f)$ has sensitive dependence on initial conditions. Then, we show (in Theorem \ref{kvocki1}), that the set $\mathcal P(f)$ of periodic points in $(X,f)$ is dense in $X$ if and only if the set $\mathcal P(f^{\star})$ of periodic points in the quotient $(X/_{\sim},f^{\star})$ is dense in $X/_{\sim}$. We also show (in Theorem \ref{kvocientek}), that for special kinds of equivalence relations $\sim $ on $X$, the dynamical system $(X,f)$ is transitive if and only if the dynamical system  $(X/_{\sim},f^{\star})$ is transitive.

In the second part of the paper, we use these results to study chaotic functions on the Cantor fan. Explicitly, we show that
\begin{enumerate}
\item  there are continuous functions $f,h:C\rightarrow C$ on the Cantor fan $C$ such that 
\begin{enumerate}
     \item $h$ is a homeomorphism and $f$ is not, 
	\item  $(C,f)$ and $(C,h)$ are both chaotic in the sense of  Devaney,   
\end{enumerate}
\item  there are continuous functions $f,h:C\rightarrow C$ on the Cantor fan $C$ such that 
\begin{enumerate}
     \item $h$ is a homeomorphism and $f$ is not, 
	\item  $(C,f)$ and $(C,h)$ are both chaotic in the sense of Robinson but not in the sense of Devaney, and
\end{enumerate}
\item  there are continuous functions $f,h:C\rightarrow C$ on the Cantor fan $C$ such that 
\begin{enumerate}
     \item $h$ is a homeomorphism and $f$ is not, 
	\item  $(C,f)$ and $(C,h)$ are both chaotic in the sense of Knudsen but not in the sense of Devaney.   
\end{enumerate}
\end{enumerate}
 In addition, we show in the third part of the paper that there are continuous functions $f,h:L\rightarrow L$ on the Lelek fan $L$ such that
\begin{enumerate}
     \item[(a)] $h$ is a homeomorphism and $f$ is not, 
	\item[(b)]  $(L,f)$ and $(L,h)$ are both chaotic in the sense of Robinson but not in the sense of Devaney.
\end{enumerate}

%It is a well-known fact that for a continuous function $f:[0,1]\rightarrow [0,1]$, the transitivity of $([0,1],f)$ implies that the set $\mathcal P(f)$ of periodic points in $([0,1],f)$ is dense in $[0,1]$, see \cite{berglund} and \cite{crannel} for more information. However, there are examples of continua $X$ and functions $f:X\rightarrow X$ such that $(X,f)$ is transitive but the set $\mathcal P(f)$ of periodic points in $(X,f)$ is not dense in $X$. Even more, there are examples of continua $X$ and functions $f:X\rightarrow X$ such that $(X,f)$ is transitive but it does not have sensitive dependence on initial conditions and the set $\mathcal P(f)$ of periodic points in $(X,f)$ is not dense in $X$. For example, if $f:S^1\rightarrow S^1$ is an irrational rotation on the circle $S^1$, then $(S^1,f)$ is an example of such a dynamical system. Here, we show that there is a homeomorphism $f:L\rightarrow L$ on the Lelek fan $L$ with top $v$ such that 
%	\begin{enumerate}
%		\item $(L,f)$ is transitive,
%		\item $(L,f)$ has sensitive dependence on initial conditions, and
%		\item $\mathcal P(f)=\{v\}$ and, therefore, the set $\mathcal P(f)$ of periodic points in $(L,f)$ is not dense in $L$.
%	\end{enumerate}
%We also show that there is a homeomorphism $g:C\rightarrow C$ on the Cantor fan $C$ with top $v$ such that 
%	\begin{enumerate}
%		\item $(C,g)$ is transitive,
%		\item $(C,g)$ has sensitive dependence on initial conditions, and
%		\item $\mathcal P(g)=\{v\}$ and, therefore, the set $\mathcal P(g)$ of periodic points in $(C,g)$ is not dense in $C$.
%	\end{enumerate}

We proceed as follows. In Section \ref{s1}, the basic definitions and results that are needed later in the paper are presented. In Section \ref{s2}, we formulate and prove Theorems \ref{kvocki}, \ref{kvocki1} and \ref{kvocientek}.  In Section \ref{s3}, we study chaos on the Cantor fan and in Section \ref{s4}, we study chaos on the Lelek fan and give some open problems.
\section{Definitions and Notation}\label{s1}
The following definitions, notation and well-known results are needed in the paper.

%\begin{definition}
%Let $X$ and $Y$ be metric spaces, and let $f:X\rightarrow Y$ be a function.  We use  $\Gamma(f)=\{(x,y)\in X\times Y \ | \ y=f(x)\}$
%to denote \emph{  the graph of the function $f$}.
%\end{definition}

\begin{definition}
Let $(X,d)$ be a metric space, $x\in X$ and $\varepsilon>0$. We use $B_d(x,\varepsilon)$ or just $B(x,\varepsilon)$ to denote the open ball,  {centered} at $x$ with radius $\varepsilon$. 
\end{definition}
\begin{definition}
We use $\mathbb N$ to denote the set of positive integers and we use $\mathbb Z$ to denote the set of integers.  
\end{definition}
%\begin{definition}
%Let $(X,d)$ be a compact metric space. Then we define \emph{$2^X$} by 
%$$
%2^{X}=\{A\subseteq X \ | \ A \textup{ is a non-empty closed subset of } X\}.
%$$
%Let $\varepsilon >0$ and let $A\in 2^X$. Then we define  \emph{$N_d(\varepsilon,A)$} by 
%$$
%N_d(\varepsilon,A)=\bigcup_{a\in A}B(a,\varepsilon).
%$$
%Let $A,B\in 2^X$. The function \emph{$H_d:2^X\times 2^X\rightarrow \mathbb R$}, defined by
%$$
%H_d(A,B)=\inf\{\varepsilon>0 \ | \ A\subseteq N_d(\varepsilon,B), B\subseteq N_d(\varepsilon,A)\},
%$$
%is called \emph{the Hausdorff metric}. The Hausdorff metric is in fact a metric and the metric space $(2^X,H_d)$ is called \emph{the hyperspace of the space $(X,d)$}. 
%\end{definition}
%\begin{remark}
%Let $(X,d)$ be a compact metric space, let $A$ be a non-empty closed subset of $X$,  and let $(A_n)$ be a sequence of non-empty closed subsets of $X$. When we say $\displaystyle A=\lim_{n\to \infty}A_n$ with respect to the {Hausdorff} metric, we mean $\displaystyle A=\lim_{n\to \infty}A_n$ in $(2^X,H_d)$. 
%\end{remark}
\begin{definition}
 \emph{A continuum} is a non-empty compact connected metric space.  \emph{A subcontinuum} is a subspace of a continuum, which is itself a continuum.
 \end{definition}
 
\begin{definition}
Let $X$ be a continuum.  We say that $X$ is \emph{a Cantor fan}, if $X$ is homeomorphic to the continuum $\bigcup_{c\in C}A_c$, where $C\subseteq [0,1]$ is a Cantor set and for each $c\in C$, $A_c$ is the  {convex} segment in the plane from $(0,0)$ to $(c,-1)$ (see Figure \ref{fig000}).
\begin{figure}[h!]
	\centering
		\includegraphics[width=30em]{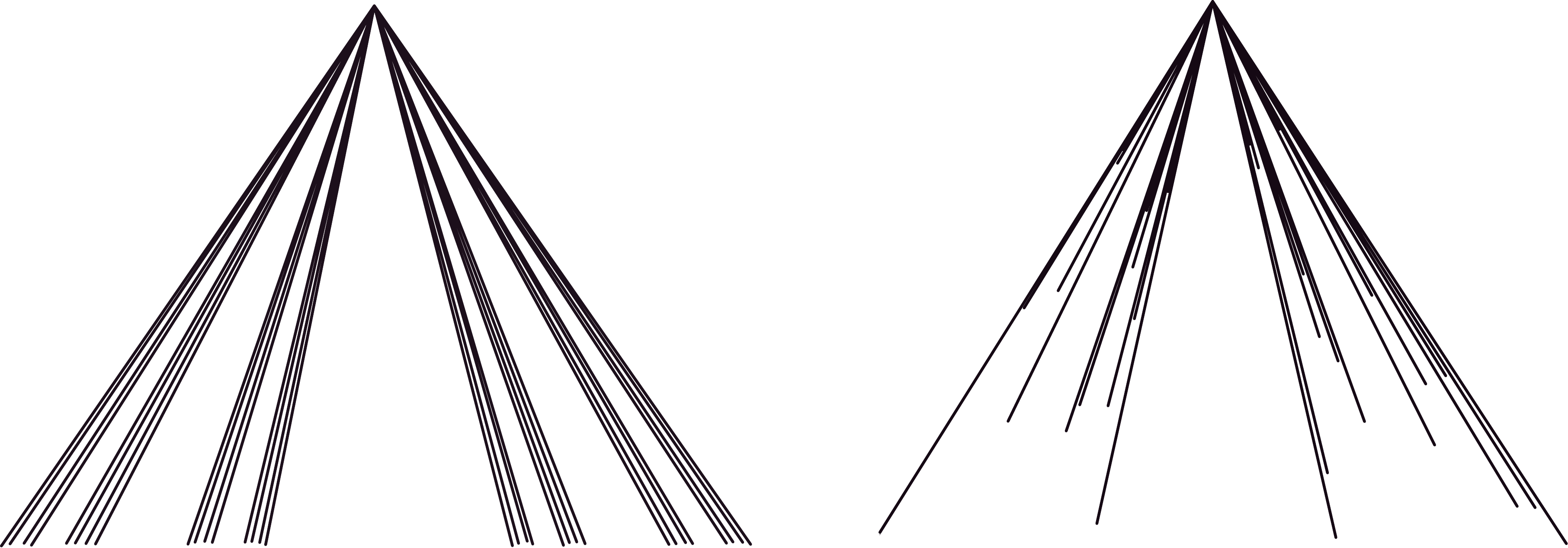}
	\caption{The Cantor fan on the left and the Lelek fan on the right}
	\label{fig000}
\end{figure}  
\end{definition}
\begin{definition}
Let $X$ be a Cantor fan and let $Y$ be a subcontinuum of $X$.  
A point $x\in Y$ is called an \emph{end-point of the continuum $Y$}, if for  every arc $A$ in $Y$ that contains $x$, $x$ is an end-point of $A$.  The set of all end-points of $Y$ will be denoted by $E(Y)$.
\end{definition}
\begin{definition} 
Let $X$ be a Cantor fan and let $Y$ be a subcontinuum of $X$.  We say that $Y$ is \emph{a Lelek fan}, if $\Cl(E(Y))=Y$ (see Figure \ref{fig000}).
\end{definition}
Note that a Lelek fan was constructed by A.~Lelek in \cite{lelek}.  An interesting  property of the Lelek fan $L$ is the fact that the set of its end-points is a dense one-dimensional set in $L$. It is also unique, i.e., any two Lelek fans are homeomorphic, for the proofs see \cite{oversteegen} and \cite{charatonik}.

b
  
\begin{definition}
Let $X$ be a {non-empty} compact metric space and let ${F}\subseteq X\times X$ be a relation on $X$. If ${F}$ is closed in $X\times X$, then we say that ${F}$ is  \emph{  a closed relation on $X$}.  
\end{definition}
\begin{definition}
Let $X$ be a {non-empty} compact metric space and let ${F}$ be a closed relation on $X$. For each positive integer $m$, we call 
$$
X_F^m=\Big\{(x_1,x_2,\ldots ,x_{m+1})\in \prod_{{ i={1}}}^{{m{ +1}}}X \ | \ \textup{ for each } i\in{  \{{1,2},\ldots ,m\}}, (x_{i},x_{i+1})\in {F}\Big\}
$$
 \emph{ the $m$-th Mahavier product of ${F}$}, we call
$$
X_F^+=\Big\{(x_1,x_2,x_3,\ldots )\in \prod_{{ i={1}}}^{\infty}X \ | \ \textup{ for each { positive} integer } i, (x_{i},x_{i+1})\in {F}\Big\}
$$
\emph{ the  Mahavier product of ${F}$}, and we call
$$
X_F=\Big\{(\ldots,x_{-3},x_{-2},x_{-1},{x_0}{ ;}x_1,x_2,x_3,\ldots )\in \prod_{i={-\infty}}^{\infty}X \ | \ \textup{ for each  integer } i, (x_{i},x_{i+1})\in {F}\Big\}
$$
\emph{the two-sided  Mahavier product of ${F}$}.
\end{definition}

\begin{definition}\label{shit}
Let $X$ be a {non-empty} compact metric space and let ${F}$ be a closed relation on $X$. 
The function  $\sigma_F^{+} : {X_F^+} \rightarrow {X_F^+}$, 
 defined by 
$$
\sigma_F^{+} ({x_1,x_2,x_3,x_4},\ldots)=({x_2,x_3,x_4},\ldots)
$$
for each $({x_1,x_2,x_3,x_4},\ldots)\in {X_F^+}$, 
is called \emph{   the shift map on ${X_F^+}$}. The function  $\sigma_F : {X_F} \rightarrow {X_F}$, 
 defined by 
$$
\sigma_F (\ldots,x_{-3},x_{-2},x_{-1},{x_0};x_1,x_2,x_3,\ldots )=(\ldots,x_{-2},x_{-1},x_{0},{x_1};x_2,x_3,x_4,\ldots )
$$
for each $(\ldots,x_{-3},x_{-2},x_{-1},{x_0};x_1,x_2,x_3,\ldots )\in {X_F}$, 
is called \emph{   the shift map on ${X_F}$}.    
\end{definition}
\begin{observation}\label{juju}
Note that $\sigma_F$ is always a homeomorphism while $\sigma_F^+$ may not be a homeomorphism.
\end{observation}
\begin{definition}
	Let $X$ be a compact metric space and let $F$ be a closed relation on $X$. The dynamical system 
	\begin{enumerate}
		\item $(X_F^{+},\sigma_F^+)$ is called \emph{a Mahavier dynamical system}.
		\item $(X_F,\sigma_F)$ is called \emph{a two-sided Mahavier dynamical system}.
	\end{enumerate}
\end{definition}
\begin{definition}
	Let $X$ be a compact metric space and let $f:X\rightarrow X$ be a continuous function. The \emph{inverse limit} generated by $(X,f)$ is the subspace
\begin{equation*}
 \varprojlim(X,f)=\Big\{(x_{1},x_{2},x_{3},\ldots ) \in \prod_{i=1}^{\infty} X \ | \ 
\text{ for each positive integer } i,x_{i}= f(x_{i+1})\Big\}
\end{equation*}
of the topological product $\prod_{i=1}^{\infty} X$.  { The function  $\sigma : \varprojlim(X,f) \rightarrow \varprojlim(X,f)$, 
 defined by 
$$
\sigma (x_1,x_2,x_3,x_4,\ldots )=(x_2,x_3,x_4,\ldots )
$$
for each $(x_1,x_2,x_3,\ldots )\in \varprojlim(X,f)$, 
is called \emph{   the shift map on $\varprojlim(X,f)$}.    }
\end{definition}
\begin{observation}
	Note that the shift map $\sigma$ on the inverse limit $\varprojlim(X,f)$ is always a homeomorphism. Also, note that for each $(x_1,x_2,x_3,\ldots )\in \varprojlim(X,f)$,
	$$
	\sigma^{-1} (x_1,x_2,x_3,\ldots )=(f(x_1),x_1,x_2,x_3,\ldots ).
	$$
\end{observation}
Theorem \ref{povezava} gives a connection between two-sided Mahavier products $X_F$ and inverse limits $\varprojlim(X_F^{+},\sigma_F^+)$.
\begin{definition}
	Let $(X,f)$  and $(Y,g)$ be dynamical systems. We say that 
	\begin{enumerate}
		\item $(Y,g)$ is topologically conjugate to $(X,f)$, if there is a homeomorphism $\varphi:X\rightarrow Y$ such that $\varphi\circ f=g\circ \varphi$.
		\item $(Y,g)$ is topologically semi-conjugate to $(X,f)$, if there is a continuous surjection $\alpha:X\rightarrow Y$ such that $\alpha\circ f=g\circ \alpha$. 
	\end{enumerate} 
\end{definition}
\begin{theorem}\label{povezava}
Let $X$ be a compact metric space and let $F$ be a closed relation on $X$. Then the following hold.
\begin{enumerate}
	\item $\varprojlim(X_F^{+},\sigma_F^+)$ is homeomorphic to $X_F$. 
	\item $(X_F,\sigma_F^{-1})$ is topologically conjugate to $(\varprojlim(X_F^{+},\sigma_F^+),\sigma)$.
\end{enumerate}  
\end{theorem}
\begin{proof}
See \cite[Theorem 4.1]{BE}.	
\end{proof}

\begin{definition}
	Let $X$ be a compact metric space and let $f:X\rightarrow X$ be a continuous function. 
We say that $(X,f)$ is a dynamical system. 
	\end{definition}
	\begin{observation}
		Note that a dynamical system $(X,f)$ can be defined in a more general setting, i.e., $X$ does not need to be a compact metric space. However, in our paper, we always deal with compact metric spaces.  
	\end{observation}
	
	We  use special kind of projections that are defined in the following definition. 
\begin{definition}
	For each (positive) integer $i$ and for each $\mathbf x=(x_1,x_2,x_3,\ldots)\in \prod_{k=1}^{\infty}X$ (or $\mathbf x=(\ldots,x_{-2},x_{-1},x_0,x_1,x_2,\ldots)\in \prod_{k=-\infty}^{\infty}X$ or $\mathbf x=(x_1,x_2,x_3,\ldots,x_{m})\in \prod_{k=1}^{m}X$), we use $\pi_i(\mathbf x)$ or $\mathbf x(i)$ or $\mathbf x_i$ to denote the $i$-th coordinate $x_i$ of the point $\mathbf x$.
%For all non-negative integers $k$ and $\ell$ such that $k\leq \ell$,  we use $[k,\ell]$ to denote the set $\{k, k+1,k+2,\ldots,\ell\}$ and $\pi_{[k,\ell]}:\prod_{i=0}^{\infty}X\rightarrow \prod_{i=k}^{\ell}X$ to denote the projection that is defined by
%$
%\pi_{[k,\ell]}(x_0,x_1,x_2,\ldots,x_k,x_{k+1},\ldots, x_{\ell},x_{\ell +1},\ldots)=(x_k,x_{k+1},x_{k+2},\ldots, x_{\ell}).
%$
%For $k=\ell$, we use $\pi_k$ to denote the projection  $\pi_{[k,k]}$.

We also use $p_1:X\times X\rightarrow X$ and $p_2:X\times X\rightarrow X$ to denote \emph{the standard projections} defined by $p_1(s,t)=s$ and $p_2(s,t)=t$ for all $(s,t)\in X\times X$.
\end{definition}
	\begin{definition}
Let $X$ be a compact metric space, let $n$ and $m$ be positive integers, $m>1$, and let $\mathbf x\in \prod_{k=1}^{n}X$,  $\mathbf y\in \prod_{k=1}^{m}X$.  If $\mathbf x(n)=\mathbf y(1)$, then we define $\mathbf x\star \mathbf y\in \prod_{k=1}^{m+n-1}X$ as follows: 
$$
\mathbf x\star \mathbf y=\Big(\mathbf x(1),\mathbf x(2),\mathbf x(3),\ldots,\mathbf x(n),\mathbf y(2),\mathbf y(2),\mathbf y(3),\ldots,\mathbf y(m)\Big).
$$
\end{definition}
\begin{definition}
		Let $(X,f)$ be a dynamical system and $p\in X$. We say that $p$ is a periodic point in $(X,f)$ if there is a positive integer $n$ such that $f^n(p)=p$. We use $\mathcal P(f)$ to denote the set of periodic points in $(X,f)$. 
	\end{definition}
	We use the following theorem about dense sets of periodic points in $(X_F^+,\sigma_F^+)$.
		\begin{theorem}\label{masten}
		Let $X$ be a compact metric space and let $F$ be a closed relation on $X$. If for each $(x,y)\in F$, there is a positive integer $n$ and a point $\mathbf z\in X_F^n$ such that $\mathbf z(1)=y$ and $\mathbf z(n+1)=x$, then the set of periodic points $\mathcal P(\sigma_F^+)$ is dense in $X_F^+$.
	\end{theorem}
	\begin{proof}
		Let $U$ be open in $\prod_{k=1}^{\infty}X$ such that $U\cap X_F^{+}\neq \emptyset$ and let $\mathbf x\in U\cap X_F^{+}$. Also, let $k$ be a positive integer and let $U_1$, $U_2$, $U_3$, $\ldots$, $U_k$ be open sets in $X$ such that 
	$$
	\mathbf x\in U_1\times U_2\times U_3\times \ldots \times U_k\times \prod_{i=k+1}^{\infty}X\subseteq U.
	$$
 Note that for each $i\in \{1,2,3,\ldots, k\}$, $(\mathbf x(i),\mathbf x(i+1))\in F$. For each $i\in \{1,2,3,\ldots, k\}$, let $n_i$  be a positive integer and let $\mathbf z_i\in X_F^{n_i}$ be such that $\mathbf z_i(1)=\mathbf x(i+1)$ and $\mathbf z_i(n_i+1)=\mathbf x(i)$. Then, let $\mathbf q\in X_F^{n_1+n_2+n_3+\ldots+n_k+k}$ be defined as 
 \begin{align*}
 	\mathbf q=\Big(\mathbf x(1),\mathbf x(2),\mathbf x(3),\ldots,\mathbf x(k),\mathbf x(k+1)\Big) \star \mathbf z_k \star \mathbf z_{k-1}\star \mathbf z_{k-2}\star \mathbf z_{k-3}\ldots \mathbf z_{3}\star \mathbf z_{2}\star \mathbf z_{1}\Big).
 \end{align*}
 Note that $\mathbf q(1)=\mathbf q(n_1+n_2+n_3+\ldots+n_k+k+1)=\mathbf x(1)$. Finally, let  
 $$
 \mathbf p=\mathbf q\star \mathbf q\star\mathbf q\star\ldots
 $$ 
 Then $\mathbf p\in U\cap \mathcal P(\sigma_F^+)$. This proves that the set of periodic points $\mathcal P(\sigma_F^+)$ is dense in $X_F^+$. 
	\end{proof}
	The following well-known theorem that follows from \cite[Lemma 2.1]{li} shows how the existence of a dense set of periodic points is carried from $(X,f)$ to the dynamical system $(\varprojlim(X,f),\sigma^{-1})$. Since the prove is short, we include it into this paper.
	\begin{theorem}\label{ingram1}
		Let $(X,f)$ be a dynamical system and let $\sigma$ be the shift homeomorphism on $\varprojlim(X,f)$. The following statements are equivalent.
		\begin{enumerate}
			\item\label{at} The set $\mathcal P(f)$ of periodic points in $(X,f)$ is dense in $X$.
			\item\label{bt} The set $\mathcal P(\sigma^{-1})$ of periodic points in $(\varprojlim(X,f),\sigma^{-1})$ is dense in $\varprojlim(X,f)$.
		\end{enumerate} 
	\end{theorem}
	\begin{proof}
		To prove the  implication from \ref{at} to \ref{bt}, suppose that  the set $\mathcal P(f)$ of periodic points in $(X,f)$ is dense in $X$. Then $\Cl(\mathcal P(f))=X$. By \cite[Lemma 2.1]{li}, 
		$$
		\Cl(\mathcal P(\sigma^{-1}))=\varprojlim(\Cl(\mathcal P(f)),f|_{\Cl(\mathcal P(f))}).
		$$ 
		It follows that 
		$$
		\Cl(\mathcal P(\sigma^{-1}))=\varprojlim(X,f|_{X})
		$$
		 and, therefore, 
		 $$
		 \Cl(\mathcal P(\sigma^{-1}))=\varprojlim(X,f).
		 $$
		  It follows that the set $\mathcal P(\sigma^{-1})$ of periodic points in $(\varprojlim(X,f),\sigma^{-1})$ is dense in $\varprojlim(X,f)$.
		
		 To prove the implication from \ref{bt} to \ref{at}, suppose that the set $\mathcal P(\sigma^{-1})$ of periodic points in $(\varprojlim(X,f),\sigma^{-1})$ is dense in $\varprojlim(X,f)$. It follows that 
		 $$
		 \Cl(\mathcal P(\sigma^{-1}))=\varprojlim(X,f).
		 $$
		 By \cite[Lemma 2.1]{li}, 
		 $$
		 \Cl(\mathcal P(\sigma^{-1}))=\varprojlim(\Cl(\mathcal P(f)),f|_{\Cl(\mathcal P(f))}).
		 $$
		  Therefore, 
		 $$
		 \varprojlim(X,f)=\varprojlim(\Cl(\mathcal P(f)),f|_{\Cl(\mathcal P(f))}).
		 $$
		  It follows from \cite[Lemma 2.6]{nadler} that $\Cl(\mathcal P(f))=X$.
	\end{proof}

	\begin{definition}
	Let $(X,f)$ be a dynamical system. We say that $(X,f)$ is 
		\begin{enumerate}
		\item \emph{transitive}, if for all non-empty open sets $U$ and $V$ in $X$,  there is a non-negative integer $n$ such that $f^n(U)\cap V\neq \emptyset$.
\item \emph{dense orbit transitive},   if there is a point $x\in X$ such that its orbit $$\{x,f(x),f^2(x),f^3(x),\ldots\}$$ is dense in $X$.
		\end{enumerate}
	\end{definition}
	\begin{observation}\label{isolatedpoints} It is a well-known fact that if $X$ has no isolated points, then $(X,f)$ is transitive if and only if $(X,f)$ is dense orbit transitive.  See \cite{A,KS} for more information about transitive dynamical systems.
\end{observation}

\begin{observation}\label{semi}
	Let $(X,f)$  and $(Y,g)$ be dynamical systems. Note that if $(X,f)$ is transitive  and if $(Y,g)$ is topologically semi-conjugate to $(X,f)$, then also $(Y,g)$ is transitive.
\end{observation}

\begin{observation}\label{semiM}
	Let $X$ be a compact metric space and let $F$ be a closed relation on $X$ such that $p_1(F)=p_2(F)=X$. Note that  $(X_F^{+},\sigma_F^+)$ is semi-conjugate to $(X_F,\sigma_F)$: for $\alpha:X_F\rightarrow X_F^{+}$, $\alpha(\mathbf x)=(\mathbf x(1),\mathbf x(2), \mathbf x(3),\ldots)$ for any $\mathbf x\in X_F$, $\alpha \circ \sigma_F=\sigma_{F}^+\circ \alpha$. 
\end{observation}
The following theorem is a well-known result; since the proof is short, we give it here for the completeness of the paper.
\begin{theorem}\label{ingram2}
		Let $(X,f)$ be a dynamical system and let $\sigma$ be the shift homeomorphism on $\varprojlim(X,f)$. If $f$ is surjective, then the following statements are equivalent.
		\begin{enumerate}
			\item\label{ats} The dynamical system $(X,f)$ is transitive.
			\item\label{bts} The dynamical system $(\varprojlim(X,f),\sigma^{-1})$ is transitive.
		\end{enumerate} 
	\end{theorem}
	\begin{proof}
	Suppose that $f$ is surjective. 		The implication from \ref{ats} to \ref{bts} follows from \cite[Theorem 2.15]{banic2}. To prove the implication from \ref{bts} to \ref{ats}, suppose that $(\varprojlim(X,f),\sigma^{-1})$ is transitive. To see that $(X,f)$ is transitive, let $U$ and $V$ be any non-empty open sets in $X$. Also, let $\hat U=U\times \prod_{k=2}^{\infty}X$ and $ \hat V=V\times \prod_{k=2}^{\infty}X$. Then $\hat U$ and $\hat V$ are non-empty open sets in $(\prod_{k=1}^{\infty}X)$ such that $\hat U\cap \varprojlim(X,f)\neq \emptyset$ and $\hat V\cap \varprojlim(X,f)\neq \emptyset$. Let $n$ be a positive integer such that $(\sigma^{-1})^n(\hat U)\cap \hat V\neq \emptyset$. Such a positive integer $n$ does exist since $(\varprojlim(X,f),\sigma^{-1})$ is transitive. Let $\mathbf x=(x_1,x_2,x_3,\ldots)\in (\sigma^{-1})^n(\hat U)\cap \hat V$. Then $x_1\in f^{n}(U)\cap V$ and, therefore,  $f^{n}(U)\cap V\neq \emptyset$. It follows that $(X,f)$ is transitive. 
			\end{proof}

	We dedicate the last part of this section to sensitive dependence on initial conditions for dynamical systems. 
	\begin{definition}
		Let $(X,f)$ be a dynamical system. We say that $(X,f)$ has \emph{sensitive dependence on initial conditions}, if there is an $\varepsilon>0$ such that for each $x\in X$ and for each $\delta >0$, there are $y\in B(x,\delta)$ and a positive integer $n$ such that  
		$$
		d(f^n(x),f^n(y))>\varepsilon. 
		$$ 
			\end{definition}
			Theorem \ref{afna} is a well-known result but we give the proof since it is not complicated and we use the results later in the paper. It enables us to work with sensitive dependence on initial conditions in a more topological setting. 			
			\begin{theorem}\label{afna}
				Let $(X,f)$ be a dynamical system. The following statements are equivalent.
				\begin{enumerate}
					\item \label{1a} $(X,f)$ has sensitive dependence on initial conditions.
					\item \label{1b} There is $\varepsilon >0$ such that for each non-empty open set $U$ in $X$, there is a positive integer $n$ such that $\diam (f^n(U))>\varepsilon$.
					\item \label{1c} There is $\varepsilon >0$ and there is a countable base $\mathcal B=\{B_1,B_2,B_3,\ldots\}$ of the space $X$ such that for each positive integer $k$, there is a positive integer $n$ such that $\diam (f^n(B_k))>\varepsilon$.
					\item \label{1d} There is $\varepsilon >0$ such that for each countable base $\mathcal B=\{B_1,B_2,B_3,\ldots\}$ of the space $X$, it holds that for each positive integer $k$, there is a positive integer $n$ such that $\diam (f^n(B_k))>\varepsilon$.
				\end{enumerate}
			\end{theorem}
			\begin{proof}
				To prove the implication from \ref{1a} to \ref{1b}, suppose that $(X,f)$ has sensitive dependence on initial conditions and let $\varepsilon>0$ be such that for each $x\in X$ and for each $\delta >0$, there are $y\in B(x,\delta)$ and a positive integer $n$ such that  $d(f^n(x),f^n(y))>\varepsilon$. We show that for each non-empty open set $U$ in $X$, there is a positive integer $n$ such that $\diam (f^n(U))>\varepsilon$.
Let $U$ be any non-empty open set in $X$, let $x\in U$ and let $\delta>0$ be such that $B(x,\delta)\subseteq U$.  It follows from our assumption that   there are $y\in B(x,\delta)$ and a positive integer $n$ such that  $d(f^n(x),f^n(y))>\varepsilon$. Since $x,y\in U$, it follows that
$$
\diam (f^n(U))\geq d(f^n(x),f^n(y))>\varepsilon
$$
and we are done. 

To prove the implication from \ref{1b} to \ref{1a}, let $\varepsilon_0 >0$ be such that for each non-empty open set $U$ in $X$, there is a positive integer $n$ such that
	$\diam (f^n(U))>\varepsilon_0$. Also, let $\displaystyle\varepsilon=\frac{\varepsilon_0}{2}$. We show that for each $x\in X$ and for each $\delta >0$, there are $y\in B(x,\delta)$ and a positive integer $n$ such that  $	d(f^n(x),f^n(y))>\varepsilon$. Let $x\in X$, let $\delta >0$, let $U=B(x,\delta)$, and let $n$ be a positive integer such that 
$\diam (f^n(U))>\varepsilon_0$. Also, let $z_1,z_2\in U$ be such that $d(f^n(z_1),f^n(z_2))>\varepsilon_0$. We consider the following possible cases.
\begin{enumerate}
	\item $z_1=x$ or $z_2=x$. Then let $y\in \{z_1,z_2\}\setminus \{x\}$ and it follows that  
	    $$
		d(f^n(x),f^n(y))=d(f^n(z_1),f^n(z_2))>\varepsilon_0>\frac{\varepsilon_0}{2}=\varepsilon. 
		$$
	\item $z_1\neq x$ and $z_2\neq x$. Suppose that $d(f^n(z_1),f^n(x))<\varepsilon$ and $d(f^n(z_2),f^n(x))<\varepsilon$. Then 
	$$
	d(f^n(z_1),f^n(z_2))\leq d(f^n(z_2),f^n(x))+d(f^n(x),f^n(z_2)) <\varepsilon+\varepsilon=2\varepsilon=\varepsilon_0,
	$$ 
	which is a contradiction. Therefore, $d(f^n(z_1),f^n(x))\geq \varepsilon$ or $d(f^n(z_2),f^n(x))\geq\varepsilon$. Suppose that  $d(f^n(z_1),f^n(x))\geq \varepsilon$. Then let $y=z_1$ and we are done.  
\end{enumerate}

To prove the implication from \ref{1b} to \ref{1c}, let $\varepsilon >0$ be such that for each non-empty open set $U$ in $X$, there is a positive integer $n$ such that
					$\diam (f^n(U))>\varepsilon$.
Let $\mathcal B=\{B_1,B_2,B_3,\ldots\}$ be any countable base of the space $X$. Since for each positive integer $k$, $B_k$ is a non-empty open set in $X$, it follows from our assumption that for each positive integer $k$, there is a positive integer $n$ such that $\diam (f^n(B_k))>\varepsilon$. 

To prove the implication from \ref{1c} to \ref{1b}, let $\varepsilon >0$ and let $\mathcal B=\{B_1,B_2,B_3,\ldots\}$ be a countable base of the space $X$ such that for each positive integer $k$, there is a positive integer $n$ such that $\diam (f^n(B_k))>\varepsilon$. We show that for each non-empty open set $U$ in $X$, there is a positive integer $n$ such that $\diam (f^n(U))>\varepsilon$.
	Let $U$ be any non-empty open set in $X$. Since $\mathcal B$ is a base in the space $X$, there is a positive integer $k$ such that $B_k\subseteq U$. It follows from our assumption that there is a positive integer $n$ such that $\diam (f^n(B_k))>\varepsilon$. Choose and fix such a positive integer $n$. Then $\diam (f^n(U))\geq \diam (f^n(B_k))>\varepsilon$.

The proof of the implication from \ref{1b} to \ref{1d} is equivalent to the proof of the implication from \ref{1b} to \ref{1c}. 

To prove the implication from \ref{1d} to \ref{1b}, let $\varepsilon >0$ and let $\mathcal B=\{B_1,B_2,B_3,\ldots\}$ be a countable base of the space $X$. It follows from our assumption that for each positive integer $k$, there is a positive integer $n$ such that $\diam (f^n(B_k))>\varepsilon$.
	The rest of the proof is equivalent to the proof of the implication from \ref{1d} to \ref{1c}. 
			\end{proof}
%			\begin{corollary}\label{alfa}
%				Let $(X,f)$ be a dynamical system and let $\mathcal B=\{B_1,B_2,B_3,\ldots\}$ be a countable base of the space $X$. The following statements are equivalent.
%				\begin{enumerate}
%					\item \label{11} $(X,f)$ has sensitive dependence on initial conditions.
%					\item \label{22} There is $\varepsilon >0$ such that for each positive integer $k$, there is a positive integer $n$ such that $	\diam (f^n(B_k))>\varepsilon$.
%									\end{enumerate}
%			\end{corollary}
%			\begin{proof}
%				The corollary follows from the fact that \ref{1a}  from Theorem \ref{afna} is equivalent to \ref{1c} from Theorem \ref{afna}.
%			\end{proof}
			We conclude this section by defining three different types of chaos.
	\begin{definition}
		Let $(X,f)$ be a dynamical system. We say that $(X,f)$ is chaotic in the sense of Robinson \cite{robinson}, if
		\begin{enumerate}
			\item  $(X,f)$ is transitive, and
			\item  $(X,f)$ has sensitive dependence on initial conditions.
		\end{enumerate}
	\end{definition}
	\begin{definition}
		Let $(X,f)$ be a dynamical system. We say that $(X,f)$ is chaotic in the sense of Knudsen \cite{knudsen}, if
		\begin{enumerate}
			\item  $\mathcal P(f)$ is dense in $X$, and
			\item  $(X,f)$ has sensitive dependence on initial conditions.
		\end{enumerate}
	\end{definition}
	\begin{definition}
		Let $(X,f)$ be a dynamical system. We say that $(X,f)$ is chaotic in the sense of Devaney \cite{devaney},  if
		\begin{enumerate}
			\item  $(X,f)$ is transitive, and
			\item  $\mathcal P(f)$ is dense in $X$.
		\end{enumerate}
	\end{definition}
\begin{observation}
	Note that it is proved in \cite{banks} that for any dynamical system $(X,f)$, if $(X,f)$ is transitive and if the set $\mathcal P(f)$ is dense in $X$, then $(X,f)$ has sensitive dependence on initial conditions. %Therefore, if $(X,f)$ is chaotic in the sense of Devaney, then it is chaotic in the sense of Robinson as well as in the sense of Knudsen.
\end{observation}
	\section{Quotients of dynamical systems}\label{s2}
	Recently, many interesting examples of dynamical systems have been produced in such a way that from a given dynamical system $(X,f)$, a new dynamical system was obtained by defining an equivalence relation on $X$ and transforming the function $f$ according to the defined equivalence relation;  see \cite{BE,BE1,BE2} for such examples. This is one of the reasons for a deeper study of quotients of dynamical systems. We begin this section by giving some basic properties of quotient spaces that will be used later. 
	
	\begin{definition}
Let $X$ be a compact metric space and let $\sim$ be an equivalence relation on $X$. For each $x\in X$, we use $[x]$ to denote the equivalence class of the element $x$ with respect to the relation $\sim$. We also use $X/_{\sim}$ to denote the quotient space $X/_{\sim}=\{[x] \ | \ x\in X\}$. 
\end{definition}
\begin{observation}\label{kvokvo}
	Let $X$ be a compact metric space, let $\sim$ be an equivalence relation on $X$, let $q:X\rightarrow X/_{\sim}$ be the quotient map that is defined by $q(x)=[x]$ for each $x\in X$,  and let $U\subseteq X/_{\sim}$. Then 
	$$
	U \textup{ is open in } X/_{\sim} ~~~ \Longleftrightarrow ~~~ q^{-1}(U)\textup{ is open in } X.
	$$
\end{observation}
\begin{definition}
Let $X$ be a compact metric space, let $\sim$ be an equivalence relation on $X$,  and let $f:X\rightarrow X$ be a function such that for all $x,y\in X$,
$$
x\sim y  \Longleftrightarrow f(x)\sim f(y).
$$
 Then we let $f^{\star}:X/_{\sim}\rightarrow X/_{\sim}$ be defined by   
$
f^{\star}([x])=[f(x)]
$
for any $x\in X$. 
\end{definition}
\begin{observation}\label{mutula}
	Let $(X,f)$ be a dynamical system. Note that we have defined a dynamical system as a pair of a compact metric space with a continuous function on it and that in this case, $X/_{\sim}$ is not necessarily metrizable. So, if $X/_{\sim}$ is metrizable, then also $(X/_{\sim},f^{\star})$ is a dynamical system. Note that in this case, $X/_{\sim}$ is semi-conjugate to $X$: for $\alpha:X\rightarrow X/_{\sim}$, defined by $\alpha(x)=q(x)$ for any $x\in X$, where $q$ is the quotient map obtained from $\sim$, $\alpha\circ f=f^{\star}\circ \alpha$.
\end{observation}
%The following well-known theorem gives some basic conditions for metrizability of the quotient space $X/_{\sim}$.
%	\begin{theorem}\label{tanasa}
%		Let $X$ be a compact metric space and let $\sim$ be a closed equivalence relation on $X$ such that for each $x\in X$, the equivalence class $[x]$ of $x$ is a closed subset of $X$. Then the quotient space $X/_{\sim}$ is metrizable. 
%	\end{theorem}
%	\begin{proof}
%		See \cite[Theorem 4.2.13]{engelking1}.
%	\end{proof}
\begin{definition}
	Let $(X,f)$ be a dynamical system and let $\sim$ be an equivalence relation on $X$  such that for all $x,y\in X$,
$$
x\sim y  \Longleftrightarrow f(x)\sim f(y).
$$
Then we say that $(X/_{\sim},f^{\star})$ is \emph{a quotient of the dynamical system $(X,f)$} or it is \emph{the quotient of the dynamical system $(X,f)$ that is obtained from the relation $\sim$}. 
\end{definition}
In the following example, we show that there is a dynamical system $(X,f)$ that has sensitive dependence on initial conditions such that its quotient $(X/_{\sim},f^{\star})$ does not, even though $X/_{\sim}$ is metrizable.  
	\begin{example}\label{buda}
		Let $(X,f)$ be a dynamical system which has sensitive dependence on initial conditions and let $\sim$ be the equivalence relation on $X$, defined by 
		$$
		x\sim y ~~~  \Longleftrightarrow  ~~~ x,y\in X  
		$$
		for all $x,y\in X$. Also, let $x_0\in X$. Then $X/_{\sim}=\{[x_0]\}$, therefore, $(X/_{\sim},f^{\star})$ does not have sensitive dependence on initial conditions. 
	\end{example}
	In Theorem \ref{kvocki},  we show that under some basic additional assumptions, the quotient $(X/_{\sim},f^{\star})$ always has sensitive dependence on initial conditions, if $(X,f)$ has sensitive dependence on initial conditions. We use the following lemmas in its proof.

	\begin{lemma}\label{trikotnik}
		Let $(X,d)$ be a compact metric space and let $A$ and $B$ be non-empty closed subsets of $X$ and let $c\in X$. Then $d(A,B)\leq d(A,c)+d(c,B)$\footnote{$d(A,B)=\inf\{d(a,b) \ | \ a\in A, b\in B\}$, $d(A,c)=\inf\{d(a,c) \ | \ a\in A\}$ and $d(c,B)=\inf\{d(c,b) \ | \ b\in B\}$.}.
	\end{lemma}
	\begin{proof}
		Let $a\in A$ be such that $d(a,c)= d(A,c)$ and let $b\in B$ be such that $d(c,b)= d(c,B)$. Then $d(A,B)\leq d(a,b)\leq d(a,c)+d(c,b)=d(A,c)+d(c,B)$. 
	\end{proof}
	\begin{lemma}\label{metrika}
		Let $(X,d)$ be a compact metric space, let $A$ be a closed subset of $X$, and let $\sim$ be the equivalence relation on $X$, defined by 
		$$
		x\sim y ~~~  \Longleftrightarrow  ~~~ x=y \textup{ or } x,y\in A  
		$$
		for all $x,y\in X$. Also, let $D:(X/_{\sim})\times (X/_{\sim})\rightarrow \mathbb R$ be defined by  
	$$
	D([x],[y])=\begin{cases}
				0\text{;} & x,y\in A \\
				d(x,A)\text{;} & x\not \in A, y\in A\\
				d(A,y)\text{;} & x\in A, y\not \in A\\
				\min\{d(x,y), d(x,A)+d(A,y)\}\text{;} & x\not \in A, y\not \in A
			\end{cases}
	$$
	for all $x,y\in X$. Then $D$ is a metric on $X/_{\sim}$. 
	\end{lemma}
	\begin{proof}
		First, note that it follows from the definition of $D$ that for all $x,y\in X$,
	\begin{enumerate}
		\item $D([x],[y])\geq 0$,
		\item $D([x],[y])=0$ if and only if $[x]=[y]$, and
		\item $D([y],[x])=D([x],[y])$.
	\end{enumerate}
	Let $x,y,z\in X$. To prove that $D([x],[y])\leq D([x],[z])+D([z],[y])$, we consider the following cases.
	\begin{enumerate}
	\item $D([x],[y])=0$. Then $D([x],[y])\leq D([x],[z])+D([z],[y])$ is obvious.
	\item $D([x],[z])=0$. Then $[x]=[z]$ and it follows that 
	$$
	D([x],[y])=D([z],[y]) = D([x],[z])+ D([z],[y]).
	$$ 
	\item $D([z],[y])=0$. Here, the proof is analogous to the proof for $D([x],[z])=0$.
	\item $D([x],[y])\neq 0$, $D([x],[z])\neq 0$ and $D([z],[y])\neq 0$. It follows that $x\neq y$, $x\neq z$ and $z\neq y$. Therefore, $\{x,y\}\not\subseteq A$, $\{x,z\}\not\subseteq A$,  and $\{z,y\}\not\subseteq A$. 
	We consider the following subcases.
	\begin{enumerate}
	\item $x\not\in A$, $y\not\in A$ and $z\in A$. Then 
	$$
	D([x],[z])+D([z],[y])=D(x,A)+D(A,y)\geq D([x],[y]).
	$$
	\item $x\in A$, $y\not\in A$ and $z\not \in A$. Then $D([x],[y])=d(A,y)$, $D([x],[z])=d(A,z)$ and $D([z],[y])=\min\{d(z,y), d(z,A)+d(A,y)\}$. We consider the following two subcases.
		\begin{enumerate}
		\item $d(z,y)\leq d(z,A)+d(A,y)$. Then $D([z],[y])=d(z,y)$. It follows from Lemma \ref{trikotnik} (taking $B=\{y\}$ and $c=z$) that $d(A,y)\leq d(A,z)+d(z,y)$. Therefore,  $D([x],[y])\leq D([x],[z])+D([z],[y])$.
		\item $d(z,y)\geq d(z,A)+d(A,y)$. Then $D([z],[y])=d(z,A)+d(A,y)$. Since $d(A,y)\leq d(A,z)+d(z,A)+d(A,y)$, it follows that the inequality $D([x],[y])\leq D([x],[z])+D([z],[y])$ holds.
		\end{enumerate}
  \item $x\not\in A$, $y\in A$ and $z\not \in A$. Here, the proof is analogous to the proof for $x\in A$, $y\not\in A$ and $z\not \in A$.
  \item $x\not\in A$, $y\not \in A$ and $z\not \in A$. We consider the following  subcases.
		\begin{enumerate}
		\item $d(x,z)\leq d(x,A)+d(A,z)$ and $d(y,z)\leq d(y,A)+d(A,z)$. It follows that 
		$$
		D([x],[y])\leq d(x,y) \leq d(x,z)+d(z,y) = D([x],[z])+D([z],[y]).
		$$
		\item  $d(x,z)\geq d(x,A)+d(A,z)$ and $d(y,z)\geq d(y,A)+d(A,z)$. Then 
		\begin{align*}
			&D([x],[y])=\min\{d(x,y), d(x,A)+d(A,y)\}\leq d(x,A)+d(A,y) \leq \\
			& d(x,A)+d(A,z)+d(z,A)+d(A,y)=D([x],[z])+D([z],[y]).
		\end{align*}
		\item  $d(x,z)\leq d(x,A)+d(A,z)$ and $d(y,z)\geq d(y,A)+d(A,z)$. Then it follows from Lemma \ref{trikotnik} that $d(x,A)\leq d(x,z)+d(z,A)$. Therefore,
		\begin{align*}
			&D([x],[y])=\min\{d(x,y), d(x,A)+d(A,y)\}\leq d(x,A)+d(A,y) \leq \\
			& d(x,z)+d(z,A)+d(A,y)=D([x],[z])+D([z],[y]).
		\end{align*}
		\item  $d(x,z)\geq d(x,A)+d(A,z)$ and $d(y,z)\leq d(y,A)+d(A,z)$. Here, the proof is analogous to the proof for $d(x,z)\leq d(x,A)+d(A,z)$ and $d(y,z)\geq d(y,A)+d(A,z)$. 
		\end{enumerate}
	\end{enumerate} 
	\end{enumerate}
	This proves that $D$ is a metric on $X/_{\sim}$.
	\end{proof}
	\begin{definition}
		Let $(X,d)$ be a compact metric space, let $A$ be a closed subset of $X$, and let $\sim$ be the equivalence relation on $X$, defined by 
		$$
		x\sim y ~~~  \Longleftrightarrow  ~~~ x=y \textup{ or } x,y\in A  
		$$
		for all $x,y\in X$. Then we use $D_A$ to denote the metric $D_A:(X/_{\sim})\times (X/_{\sim})\rightarrow \mathbb R$, which is defined by  
	$$
	D_A([x],[y])=\begin{cases}
				0\text{;} & x,y\in A\\
				d(x,A)\text{;} & x\not \in A, y\in A\\
				d(A,y)\text{;} & x\in A, y\not \in A\\
				\min\{d(x,y), d(x,A)+d(A,y)\}\text{;} & x\not \in A, y\not \in A
			\end{cases}
	$$
	for all $x,y\in X$. 
	\end{definition}
	\begin{observation}
		Let $(X,f)$ be a dynamical system, let $A$ be a closed subset of $X$ such that $f(A)\subseteq A$ and $f(X\setminus A)\subseteq X\setminus A$, 
		 and let $\sim$ be the equivalence relation on $X$, defined by 
		$$
		x\sim y ~~~  \Longleftrightarrow  ~~~ x=y \textup{ or } x,y\in A  
		$$
		for all $x,y\in X$. Note that for all $x,y\in X$, 
		$$
x\sim y  \Longleftrightarrow f(x)\sim f(y).
$$
	\end{observation}
	\begin{lemma}\label{topologiji}
		Let $(X,d)$ be a compact metric space, let $A$ be a closed subset of $X$, and let $\sim$ be the equivalence relation on $X$, defined by 
		$$
		x\sim y ~~~  \Longleftrightarrow  ~~~ x=y \textup{ or } x,y\in A  
		$$
		for all $x,y\in X$. Also, let
		\begin{enumerate}
			\item $\mathcal T_d$ be the topology on $X$ such that the family $\mathcal B_d=\{B_d(x,r) \ | \  x\in X, r>0\}$ of open balls $B_d(x,r)$ in $X$ is the base for $\mathcal T_d$,
			\item $\mathcal T_{D_A}$ be the topology on $X/_{\sim}$ such that the family $\mathcal B_{D_A}=\{B_{D_A}([x],r) \ | \  x\in X, r>0\}$ of open balls  $B_{D_A}([x],r)$ in $X/_{\sim}$ is the base for $\mathcal T_{D_A}$,
			\item $q:X\rightarrow X/_{\sim}$ be the quotient map, and 
			\item $\mathcal Q$ be the quotient topology on $X/_{\sim}$, obtained from topology $\mathcal T_d$; i.e.,  for each $U\subseteq X/_{\sim}$,
		$$
	U \in \mathcal Q ~~~ \Longleftrightarrow ~~~ q^{-1}(U)\in \mathcal T_d.
	$$
	\end{enumerate}
		Then $\mathcal T_{D_A}=\mathcal Q$.
	\end{lemma}
	\begin{proof}
		In the proof, we use $D$ to denote the metric $D=D_A$. To see that $\mathcal T_D\subseteq \mathcal Q$, we show that $\mathcal B_D\subseteq \mathcal Q$. Let $x\in X$ and $r>0$ be arbitrarily chosen. We show that $B_D([x],r)\in \mathcal Q$ by showing that $q^{-1}(B_D([x],r))\in \mathcal T_d$. We consider the following possible cases.
\begin{enumerate}
	\item $x\in A$. Then
	\begin{align*}
	B_D([x],r)=\{[y]\in X/_{\sim} \ | \ D([x],[y])<r\}=\{[y]\in X/_{\sim} \ | \ d(A,y)<r\}.
\end{align*}
It follows that 
\begin{align*}
	q^{-1}(B_D([x],r))=&q^{-1}(\{[y]\in X/_{\sim} \ | \ d(A,y)<r\})=\{y\in X \ | \ d(A,y)<r\}=\\
	&\bigcup_{a\in A}B_d(a,r).
\end{align*}
Therefore, $q^{-1}(B_D([x],r))\in \mathcal T_d$.
	\item $x\not \in A$ and $r\leq d(x,A)$. Note that for each $a\in A$,
	\begin{align*}
		D([a],[x])=d(A,x)\geq r
	\end{align*}
	and, therefore, $[a]\not\in B_D([x],r)$. It follows that for each $[y]\in B_D([x],r)$, $\min\{d(x,y), d(x,A)+d(A,y)\}=d(x,y)$ and, therefore,  $D([x],[y])=d(x,y)$. Hence,	
	\begin{align*}
	B_D([x],r)=&\{[y]\in X/_{\sim} \ | \ D([x],[y])<r\}=\{[y]\in X/_{\sim} \ | \ d(x,y)<r\}
\end{align*}
and it follows that  
\begin{align*}
	q^{-1}(B_D([x],r))=&q^{-1}(\{[y]\in X/_{\sim} \ | \ d(x,y)<r\})=\{y\in X \ | \ d(x,y)<r\}=B_d(x,r).
\end{align*}
Therefore, $q^{-1}(B_D([x],r))\in \mathcal T_d$.
\item $x\not \in A$ and $r > d(x,A)$. 
%Note that in this case, for each $a\in A$,
%	\begin{align*}
%		D([a],[x])=d(A,x)<r
%	\end{align*}
%	and, therefore, $[a]\in B_D([x],r)$. 
We show that in this case, 
	$$
	q^{-1}(B_D([x],r))=B_d(x,r)\cup \Big(\bigcup_{a\in A}B_d(a,r-d(x,A))\Big).
	$$ 
	To see that $q^{-1}(B_D([x],r))\subseteq B_d(x,r)\cup \Big(\bigcup_{a\in A}B_d(a,r-d(x,A))\Big)$, take any point $y\in q^{-1}(B_D([x],r))$. Then $q(y)=[y]\in B_D([x],r)$ and $D([x],[y])<r$ follows. We consider the following cases.
	\begin{enumerate}
	\item $y\in A$. Then $y\in B_d(y,r-d(x,A))\subseteq \bigcup_{a\in A}B_d(a,r-d(x,A))$ and we are done. 
	\item $y\not \in A$. Then $D([x],[y])=\min\{d(x,y),d(x,A)+d(A,y)\}$. If $d(x,y)\leq d(x,A)+d(A,y)$, then $D([x],[y])=d(x,y)$ and it follows that $d(x,y)<r$. Therefore, $y\in B_d(x,r)$. If $d(x,y)\geq d(x,A)+d(A,y)$, then $D([x],[y])=d(x,A)+d(A,y)$ and it follows that $d(x,A)+d(A,y)<r$. Therefore, $d(A,y)<r-d(x,A)$ and $y\in \bigcup_{a\in A}B_d(a,r-d(x,A))$ follows. 
	\end{enumerate} 
	It follows that $y\in B_d(x,r)\cup \Big(\bigcup_{a\in A}B_d(a,r-d(x,A))\Big)$. 
	
	To see that $B_d(x,r)\cup \Big(\bigcup_{a\in A}B_d(a,r-d(x,A))\Big)\subseteq q^{-1}(B_D([x],r))$, take any point  $y\in B_d(x,r)\cup \Big(\bigcup_{a\in A}B_d(a,r-d(x,A))\Big)$. 
	\begin{enumerate}
	\item $y\in A$. Then $D([x],[y])= d(x,A)<r$ and $y\in q^{-1}(B_D([x],r))$ follows.
	\item $y\not\in A$. First, suppose that $y\in B_d(x,r)$.  Then
	$$
	D([x],[y])= \min\{d(x,y), d(x,A)+d(A,y)\}\leq d(x,y)<r
	$$  
	and $y\in q^{-1}(B_D([x],r))$ follows. 
	
	Next, suppose that $y\in \bigcup_{a\in A}B_d(a,r-d(x,A))$. Then 
	\begin{align*}
	D([x],[y])= &\min\{d(x,y), d(x,A)+d(A,y)\}\leq \\
	&d(x,A)+d(A,y)<d(x,A)+r-d(x,A)=r	
	\end{align*} 
	and $y\in q^{-1}(B_D([x],r))$ follows.
	\end{enumerate}
	\end{enumerate}
We have just proved that $\mathcal T_D\subseteq \mathcal Q$. Next, we show that $\mathcal Q\subseteq \mathcal T_D$. Let $U\in \mathcal Q$. To see that $U\in \mathcal T_D$, let $[x]\in U$. It follows that $[x]\subseteq q^{-1}(U)$. Note that  $q^{-1}(U)$ is open in $X$ (since $U\in \mathcal Q$). We consider the following cases.
\begin{enumerate}
	\item $x\in A$. Then, since $[x]=A$, it follows that $A\subseteq q^{-1}(U)$. Since $A$ is compact, there is an $r>0$ such that 
	$$
	\bigcup_{a\in A}B_d(a,r)\subseteq q^{-1}(U).
	$$
	Choose and fix such an $r$. We show that
	$$
	B_D([x],r)\subseteq U. 
	$$
	Let $[y]\in B_D([x],r)$. Then $d(A,y)=D([x],[y])<r$, and therefore, $y\in \bigcup_{a\in A}B_d(a,r)$. Since $\bigcup_{a\in A}B_d(a,r)\subseteq q^{-1}(U)$, it follows that $q(y)=[y]\in U$. It follows that $U\in \mathcal T_D$.
	\item $x\not\in A$. Then let $r\in (0,d(x,A))$ be such that $B_d(x,r)\subseteq q^{-1}(U)$. We show that
	$$
	B_D([x],r)\subseteq U. 
	$$
	Let $[y]\in B_D([x],r)$. Then $d(x,y)=D([x],[y])<r$, and therefore, $y\in B_d(x,r)$. Since $B_d(x,r)\subseteq q^{-1}(U)$, it follows that $q(y)=[y]\in U$. 
	It follows that $U\in \mathcal T_D$.
\end{enumerate}
We have just proved that also $\mathcal T_D\subseteq \mathcal Q$. It follows that $\mathcal T_D=\mathcal Q$. 
	\end{proof}
	Our next goal is to prove Theorem \ref{kvocki} about sensitive dependence on initial conditions of quotients of dynamical systems. For that purpose, in the following definition, we introduce a property that is a stronger version of sensitive dependence on initial conditions.
	\begin{definition}
		Let $(X,f)$ be a dynamical system and let $A$ be a non-empty closed subset of $X$. We say that \emph{$(X,f)$ has sensitive dependence on initial conditions with respect to $A$}, if there is $\varepsilon>0$ such that for each non-empty open set $U$ in $X$, there are $x,y\in U$ and a positive integer $n$ such that 
		$$
		\min\{d(f^n(x),f^n(y)),d(f^n(x),A)+d(f^n(y),A)\}>\varepsilon.
		$$  
	\end{definition}
	\begin{observation}
		Let $(X,f)$ be a dynamical system and let $A$ be a non-empty closed subset of $X$. Note that it follows from Theorem \ref{afna} that if $(X,f)$ has sensitive dependence on initial conditions with respect to $A$, then $(X,f)$ has sensitive dependence on initial conditions.
			\end{observation}
	The following theorem shows how sensitive dependence on initial conditions with respect to a set is carried from a dynamical system $(X,f)$ to the dynamical system $(\varprojlim(X,f),\sigma^{-1})$.
	\begin{theorem}\label{ingram}
		Let $(X,f)$ be a dynamical system, let $A$ be a non-empty closed subset of $X$ such that $f(A)\subseteq A$, and let
		 $\sigma$ be the shift homeomorphism on $\varprojlim(X,f)$. If $f$ is surjective and if $(X,f)$ has sensitive dependence on initial conditions with respect to $A$, then $(\varprojlim(X,f),\sigma^{-1})$ has sensitive dependence on initial conditions with respect to $\varprojlim(A,f|_A)$. 
	\end{theorem}
	\begin{proof}
		Let $d$ be the metric on $X$ and suppose that $f$ is surjective. We use $\rho$ to denote the product metric on $\prod_{k=1}^{\infty}X$, which is defined by  		
		$$
		\rho((x_1,x_2,x_3,\ldots),(y_1,y_2,y_3,\ldots))=\max\Big\{\frac{d(x_k,y_k)}{2^k} \ \big| \  k \textup{ is a positive integer}\Big\}
		$$ 
		for all $(x_1,x_2,x_3,\ldots),(x_1,x_2,x_3,\ldots)\in \prod_{k=1}^{\infty}X$.
	
Let $\varepsilon>0$ be such that for each non-empty open set $U$ in $X$, there are $x,y\in U$ and a positive integer $m$ such that 
		$$
		\min\{d(f^m(x),f^m(y)),d(f^m(x),A)+d(f^m(y),A)\}>\varepsilon.
		$$
		Let $U$ be a basic set of the product topology on $\prod_{k=1}^{\infty}X$ such that $U\cap \varprojlim(X,f)\neq \emptyset$. Also, let $n$ be a positive integer and for each $i\in \{1,2,3,\ldots, n\}$, let $U_i$ be an open set in $X$ such that 
	$$
	U=U_1\times U_2\times U_3\times \ldots \times U_n\times \prod_{k=n+1}^{\infty}X. 
	$$
	Also, let $x,y\in U_1$ and let $m$ be a positive integer such that 
	$$
		\min\{d(f^m(x),f^m(y)),d(f^m(x),A)+d(f^m(y),A)\}>\varepsilon.
		$$
		Then let $\mathbf x,\mathbf y\in U\cap \varprojlim(X,f)$ be any points such that $\mathbf x(1)=x$ and $\mathbf y(1)=y$. Such points $\mathbf x$ and $\mathbf y$ do exist since $f$ is surjective. Then
		\begin{align*}
			&\rho\Big((\sigma^{-1})^m(\mathbf x),(\sigma^{-1})^m(\mathbf y)\Big)=\rho\Big((\sigma^{-m})(\mathbf x),(\sigma^{-m})(\mathbf y)\Big)=\\
			&\rho\Big((f^m(x),\ldots,f(x),x,\mathbf x(2), \mathbf x(3),\ldots ),(f^m(y),\ldots,f(y),y,\mathbf y(2), \mathbf y(3),\ldots )\Big)=\\
			&\max\Big(\Big\{\frac{d(f^k(x),f^k(y))}{2^{m-k+1}} \ | \  k\in \{0,1,2,3,\ldots, m\}\Big\}\cup \Big\{\frac{d(\mathbf x(k),\mathbf y(k))}{2^{k+m}} \ | \  k\in \{2,3,4,\ldots\}\Big\}\Big)\geq \\
			& \frac{d(f^m(x),f^m(y))}{2}\geq \frac{\min\{d(f^m(x),f^m(y)),d(f^m(x),A)+d(f^m(y),A)\}}{2}>\frac{\varepsilon}{2}.
		\end{align*}
		Next, let $\mathbf a,\mathbf b\in \varprojlim(A,f|_A)$ be such that $\rho\Big((\sigma^{-1})^m(\mathbf x),\varprojlim(A,f|_A)\Big)=\rho\Big((\sigma^{-1})^m(\mathbf x),\mathbf a\Big)$ and $\rho\Big((\sigma^{-1})^m(\mathbf y),\varprojlim(A,f|_A)\Big)=\rho\Big((\sigma^{-1})^m(\mathbf y),\mathbf b\Big)$. Then 
		\begin{align*}
			&\rho\Big((\sigma^{-1})^m(\mathbf x),\varprojlim(A,f|_A)\Big)+\rho\Big((\sigma^{-1})^m(\mathbf y),\varprojlim(A,f|_A)\Big)=\\
			&\rho\Big((\sigma^{-1})^m(\mathbf x),\mathbf a\Big)+\rho\Big((\sigma^{-1})^m(\mathbf y),\mathbf b\Big)=\\
			&\rho\Big((f^m(x),\ldots,f(x),x,\mathbf x(2), \mathbf x(3),\ldots ),\mathbf a\Big)+\rho\Big((f^m(y),\ldots,f(y),y,\mathbf y(2), \mathbf y(3),\ldots ),\mathbf b\Big)=\\
			&\rho \Big((f^m(x),\ldots,f(x),x,\mathbf x(2), \mathbf x(3),\ldots ),(\mathbf a(1),\mathbf a(2),\mathbf a(3),\ldots)\Big)+\\
			&\rho \Big((f^m(y),\ldots,f(y),y,\mathbf y(2), \mathbf y(3),\ldots ),(\mathbf b(1),\mathbf b(2),\mathbf b(3),\ldots)\Big)=\\
			&\max\Big(\Big\{ \frac{d(f^k(x),\mathbf a(m-k+1))}{2^{m-k+1}} \ | \  k\in \{0,1,2,\ldots, m\} \Big\}\cup \Big\{ \frac{d(\mathbf x(k),\mathbf a(m+k))}{2^{m+k}} \ | \  k\in \{2,3,4,\ldots \} \Big\}\Big)+\\
			&\max\Big(\Big\{ \frac{d(f^k(y),\mathbf b(m-k+1))}{2^{m-k+1}} \ | \  k\in \{0,1,2,\ldots, m\} \Big\}\cup \Big\{ \frac{d(\mathbf y(k),\mathbf b(m+k))}{2^{m+k}} \ | \  k\in \{2,3,4,\ldots \} \Big\}\Big)\geq \\
			&\frac{d(f^m(x),\mathbf a(1))}{2}+\frac{d(f^m(y),\mathbf b(1))}{2}=\frac{d(f^m(x),\mathbf a(1))+d(f^m(y),\mathbf b(1))}{2}\geq \\
			&\frac{\min\{d(f^m(x),f^m(y)),d(f^m(x),A)+d(f^m(y),A)\}}{2}>\frac{\varepsilon}{2}.
		\end{align*}
		This proves that $(\varprojlim(X,f),\sigma^{-1})$ has sensitive dependence on initial conditions with respect to $\varprojlim(A,f|_A)$.
		    	\end{proof}
	The following is our first  main result of this section. 
		\begin{theorem}\label{kvocki}
	Let $(X,f)$ be a dynamical system, let $A$ be a nowhere dense closed subset of $X$ such that $f(A)\subseteq A$ and $f(X\setminus A)\subseteq X\setminus A$, and let $\sim$ be the equivalence relation on $X$, defined by 
		$$
		x\sim y ~~~  \Longleftrightarrow  ~~~ x=y \textup{ or } x,y\in A  
		$$
		for all $x,y\in X$. 
The following statements are equivalent.
	\begin{enumerate}
		\item\label{uno1} $(X,f)$ has sensitive dependence on initial conditions with respect to $A$. 
		\item\label{due2} $(X/_{\sim},f^{\star})$ has sensitive dependence on initial conditions.
	\end{enumerate}	 
	\end{theorem}
	\begin{proof}
	First, let $q:X\rightarrow X/_{\sim}$ be the quotient map defined by $q(x)=[x]$ for any $x\in X$. 
	
	To prove the implication from \ref{uno1} to \ref{due2}, suppose that $(X,f)$ has sensitive dependence on initial conditions with respect to $A$. Let $\varepsilon>0$ be such that for each non-empty open set $U$ in $X$, there are $x,y\in U$ and a positive integer $n$ such that 
		$$
		\min\{d(f^n(x),f^n(y)),d(f^n(x),A)+d(f^n(y),A)\}>\varepsilon.
		$$  
Also,  let $U$ be any non-empty open set in $X/_{\sim}$. Then $q^{-1}(U)$ is open in $X$. Also, let $V=q^{-1}(U)\setminus A$. Since $A$ is closed in $X$, it follows that $V$ is open in $X$. Since $A$ is nowhere dense in $X$, it follows that $V$ is non-empty. Let $x,y\in V$ and let $n$ be a positive integer such that $\min\{d(f^n(x),f^n(y)),d(f^n(x),A)+d(f^n(y),A)\}>\varepsilon$. Note that since $f(X\setminus A)\subseteq X\setminus A$, it follows that $f^n(x)\not \in A$ and $f^n(y)\not\in A$. Therefore, 
        \begin{align*}
        	D_A\big((f^{\star})^n([x]),(f^{\star})^n([y])\big)=&D_A\big([f^n(x)],[f^n(y)]\big)=\\
        	&\min\{d(f^n(x),f^n(y)), d(f^n(x),A)+d(A,f^n(y))\}>\varepsilon.
        \end{align*}
        Therefore, $\diam((f^{\star})^n(U))>\varepsilon$ and it follows that $(X/_{\sim},f^{\star})$ has sensitive dependence on initial conditions. 
        
        Next, we prove the implication from \ref{due2} to \ref{uno1}. Suppose that $(X/_{\sim},f^{\star})$ has sensitive dependence on initial conditions. Let $\varepsilon>0$ be such that for any non-empty open set $U$ in $X/_{\sim}$, there is a positive integer $n$ such that 
     $\diam((f^{\star})^n(U))>\varepsilon$. To prove that $(X,f)$ has sensitive dependence on initial conditions with respect to $A$, let $W$ be any non-empty open set in $X$. Also, let $V=W\setminus A$. Since $A$ is closed in $X$, it follows that $V$ is open in $X$. Since $A$ is nowhere dense in $X$, it follows that $V$ is non-empty. Since $V\cap A=\emptyset$, it follows that $q(V)$ is non-empty and open in $X/_{\sim}$. Let $m$ be a positive integer such that $\diam ((f^{\star})^m(q(V))) > \varepsilon$ and let $x,y\in V$ be such that $D_A((f^{\star})^m([x]),(f^{\star})^m([y]))>\varepsilon$. Then $D_A([f^m(x)],[f^m(y)])>\varepsilon$.  Since $f(X\setminus A)\subseteq X\setminus A$, it follows that $f^m(x),f^m(y)\not \in A$ and, therefore, $D_A([f^m(x)],[f^m(y)])=\min\{d(f^m(x),f^m(y)), d(f^m(x),A)+d(A,f^m(y))\}>\varepsilon$.
     This proves that $(X,f)$ has sensitive dependence on initial conditions with respect to $A$.
         \end{proof}
     	The following theorem is our second main result of this section. 
		\begin{theorem}\label{kvocki1}
	Let $(X,f)$ be a dynamical system, let $A$ be a nowhere dense closed subset of $X$ such that $f(A)\subseteq A$ and $f(X\setminus A)\subseteq X\setminus A$, and let $\sim$ be the equivalence relation on $X$, defined by 
		$$
		x\sim y ~~~  \Longleftrightarrow  ~~~ x=y \textup{ or } x,y\in A  
		$$
		for all $x,y\in X$. 
		Then the following statements are equivalent.
		\begin{enumerate}
			\item\label{ata} The set $\mathcal P(f)$ of periodic points in $(X,f)$ is dense in $X$.
			\item\label{btb} The set $\mathcal P(f^{\star})$ of periodic points in the quotient $(X/_{\sim},f^{\star})$ is dense in $X/_{\sim}$.
		\end{enumerate} 
	\end{theorem}
	\begin{proof}
	First, let $q:X\rightarrow X/_{\sim}$ be the quotient map defined by $q(x)=[x]$ for any $x\in X$. 	To prove the implication from \ref{ata} to \ref{btb}, suppose that the set $\mathcal P(f)$ of periodic points in $(X,f)$ is dense in $X$. Let $U$ be any non-empty open set in $X/_{\sim}$. Then $q^{-1}(U)$ is open in $X$. Also, let $V=q^{-1}(U)\setminus A$. Since $A$ is closed in $X$, it follows that $V$ is open in $X$. Since $A$ is nowhere dense in $X$, it follows that $V$ is non-empty. Let $p\in V\cap \mathcal P(f)$ and let $n$ be a positive integer such that $f^n(p)=p$. Note that since $f(X\setminus A)\subseteq X\setminus A$, it follows that for each $k\in\{1,2,3,\ldots, n\}$, $f^k(x)\not \in A$. Therefore, 
        \begin{align*}
        	(f^{\star})^n([p])=[f^n(p)]=[p].
        \end{align*}
        Therefore, $[p]\in \mathcal P(f^{\star})\cap U$ and it follows that the set $\mathcal P(f^{\star})$ of periodic points in  $(X/_{\sim},f^{\star})$ is dense in $X/_{\sim}$. 
        
        To prove the implication from \ref{btb} to \ref{ata}, suppose that the set $\mathcal P(f^{\star})$ of periodic points in the quotient $(X/_{\sim},f^{\star})$ is dense in $X/_{\sim}$. Let $U$ be any non-empty open set in $X$. Also, let $V=U\setminus A$. Since $A$ is closed in $X$, it follows that $V$ is open in $X$. Since $A$ is nowhere dense in $X$, it follows that $V$ is non-empty. It follows from $q^{-1}(q(V))=V$ that $q(V)$ is a non-empty open set in $X/_{\sim}$, therefore, there is $[p]\in q(V)\cap \mathcal P(f^{\star})$.  Let $n$ be a positive integer such that $(f^{\star})^n([p])=[p]$. Note that $p\in X\setminus A$ and since $f(X\setminus A)\subseteq X\setminus A$, it follows that $f^n(p)\in X\setminus A$. Therefore,
        $$
        f^n(p)\in q^{-1}(q(f^n(p)))=q^{-1}([f^n(p)])=q^{-1}((f^{\star})^n([p]))=q^{-1}([p])= q^{-1}(q(p))=\{p\}
        $$ 
        and it follows that $p\in \mathcal P(f)$. Since $p\in V$ and $V\subseteq U$, it follows that $p\in  \mathcal P(f)\cap U$. Therefore, the set $\mathcal P(f)$ of periodic points in $(X,f)$ is dense in $X$.
         \end{proof}
We finish the section by stating and proving a theorem about transitivity. 
\begin{theorem}\label{kvocientek}
	Let $(X,f)$ be a dynamical system, let $A$ be a nowhere dense closed subset of $X$ such that $f(A)\subseteq A$ and $f(X\setminus A)\subseteq X\setminus A$, and let $\sim$ be the equivalence relation on $X$, defined by 
		$$
		x\sim y ~~~  \Longleftrightarrow  ~~~ x=y \textup{ or } x,y\in A  
		$$
		for all $x,y\in X$. 
		Then the following statements are equivalent.
		\begin{enumerate}
			\item\label{atas} The dynamical system $(X,f)$ is transitive.
			\item\label{btbs} The dynamical system $(X/_{\sim},f^{\star})$ is transitive.
		\end{enumerate} 
	\end{theorem}
	\begin{proof}
	The implication from \ref{atas} to \ref{btbs} follows from Theorem \ref{kvocienti}. To prove the implication from \ref{btbs} to \ref{atas}, suppose that $(X/_{\sim},f^{\star})$ is transitive. To see that $(X,f)$ is transitive, let $U$ and $V$ be any non-empty open sets in $X$. Also, let $\hat U=U\setminus A$ and $\hat V=V\setminus A$. Since $A$ is closed in $X$, it follows that $\hat U$ and $\hat V$ are both open in $X$. Since $A$ is nowhere dense in $X$, it follows that $\hat U$ and $\hat V$ are both non-empty. It follows from $q^{-1}(q(\hat U))=\hat U$ and $q^{-1}(q(\hat V))=\hat V$ that $q(\hat U)$ and $q(\hat V)$ are both non-empty open sets in $X/_{\sim}$. Since $(X/_{\sim},f^{\star})$ is transitive, there is a positive integer $n$ such that $(f^{\star})^n(q(\hat U))\cap q(\hat V)\neq \emptyset$. Note that for each $x\in X$, if $x\in X\setminus A$, then $f^n(x)\in X\setminus A$ (since $f(X\setminus A)\subseteq X\setminus A$). It follows that
	\begin{align*}
		&(f^{\star})^n(q(\hat U))\cap q(\hat V)=\{(f^{\star})^n(q(x)) \ | \ x\in \hat U\}\cap q(\hat V)=\\
		&\{(f^{\star})^n([x]) \ | \ x\in \hat U\}\cap q(\hat V)=\{[f^n(x)] \ | \ x\in \hat U\}\cap q(\hat V)=\\
		&\{q(f^n(x)) \ | \ x\in \hat U\}\cap q(\hat V)=q(f^n(\hat U))\cap q(\hat V)= q(f^n(\hat U)\cap \hat V)
	\end{align*}
	since $q$ restricted to $X\setminus A$ is bijective. It follows from this that $q(f^n(\hat U)\cap \hat V)\neq \emptyset$, therefore, $f^n(\hat U)\cap \hat V\neq \emptyset$. We have just proved that there is a positive integer $n$ such that $f^n(U)\cap V\neq \emptyset$. It follows that  $(X,f)$ is transitive. 
		         \end{proof}

In Section \ref{s3}, we demonstrate possible applications of Theorems  \ref{kvocki},  \ref{kvocki1}, and \ref{kvocientek} to obtain different variants of chaos on the Cantor fan.
   %%%%%%%%%
   %%%%%
   %%%
   %
 \section{Chaos on the Cantor fan}\label{s3}
 
          In this section, we present different types of chaos on the Cantor fan. In  Subsection \ref{32}, we study functions $f$ on the Cantor fan $C$ such that $(C,f)$ is chaotic in the sense of Devaney, in Subsection \ref{31}, we study functions $f$ on the Cantor fan $C$ such that $(C,f)$ is chaotic in the sense of Robinson but not in the sense of Devaney, and in Subsection \ref{33}, we study functions $f$ on the Cantor fan $C$ such that $(C,f)$  is chaotic in the sense of Knudsen but it is not chaotic in the sense of Devaney.  
        
   \subsection{Devaney's chaos on the Cantor fan}\label{32} 
  Here, we study functions $f$ on the Cantor fan $C$ such that $(C,f)$ is chaotic in the sense of Devaney.

  \begin{definition}\label{nall2}
In this subsection, we use $X$ to denote $X=[0,1]\cup [2,3]\cup [4,5]\cup [6,7]$, and we let $f_1,f_2,f_3:X\rightarrow X$ be the homeomorphisms from $X$ to $X$ that are defined by  
$$
f_1(x)=\begin{cases}
				x^2\text{;} & x\in [0,1]\\
				(x-2)^{\frac{1}{2}}+2\text{;} & x\in [2,3]\\
				(x-4)^3+4\text{;} & x\in [4,5]\\
				(x-6)^{\frac{1}{3}}+6\text{;} & x\in [6,7]
			\end{cases} 
$$
$$
f_2(x)=\begin{cases}
				x+2\text{;} & x\in [0,1]\cup  [4,5]\\
				x-2\text{;} & x\in [2,3]\cup [6,7]
			\end{cases} 
$$
$$
f_3(x)=\begin{cases}
				x^2\text{;} & x\in [0,1]\\
				x+2\text{;} & x\in [2,3]\\
				x-2\text{;} & x\in [4,5]\\
				(x-6)^{\frac{1}{3}}+6\text{;} & x\in [6,7]
			\end{cases}  
$$
for each $x\in X$. Then we use $F$ to denote the relation $F=\Gamma(f_1)\cup \Gamma(f_2)\cup \Gamma(f_3)$; see Figure \ref{figurca}.
\begin{figure}[h!]
	\centering
		\includegraphics[width=26em]{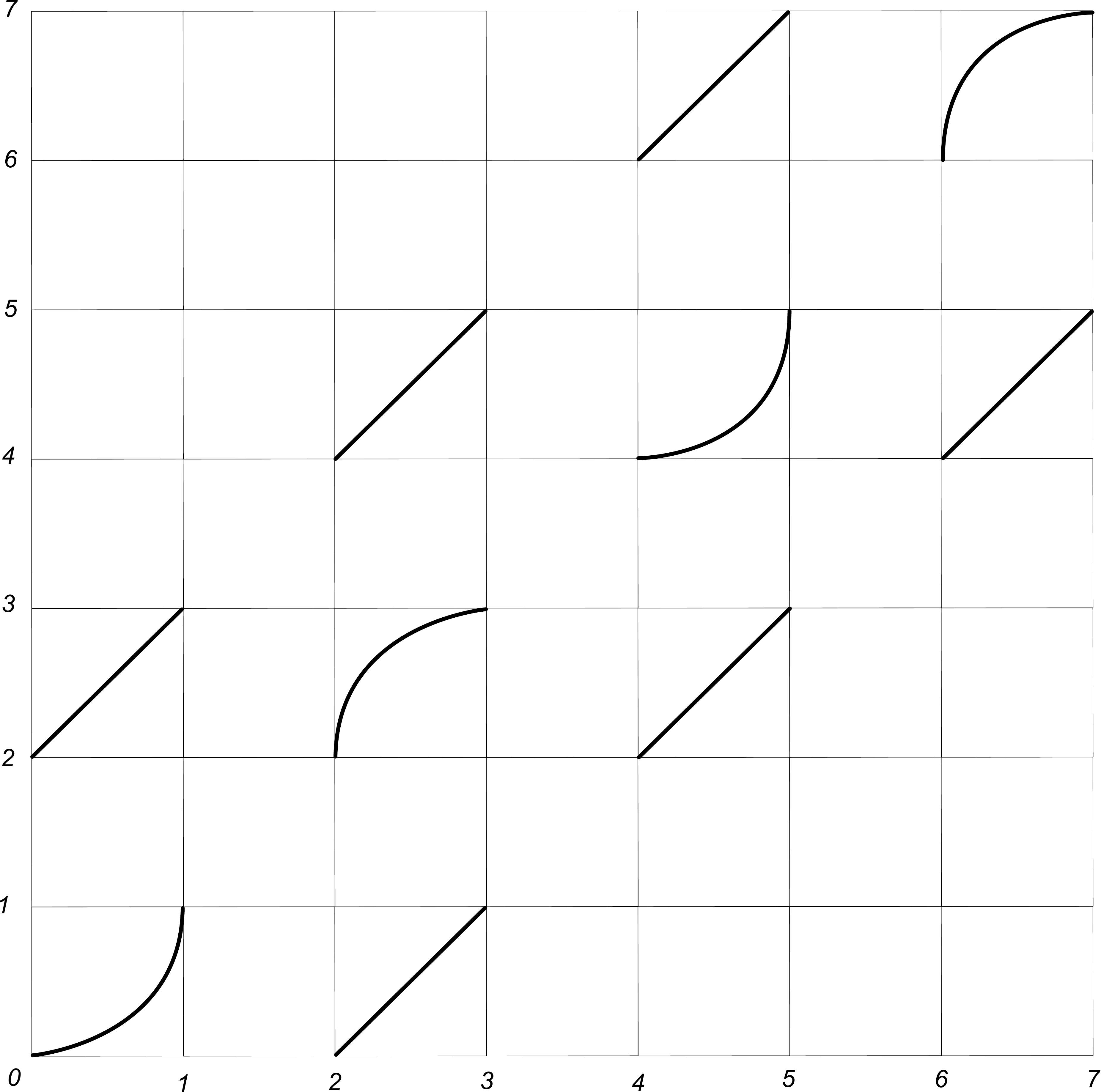}
	\caption{The relation $F$ from Definition \ref{nall2}}
	\label{figurca}
\end{figure}  
	\end{definition}
\begin{definition}
	We define two equivalence relations. 
	\begin{enumerate}
		\item For all $\mathbf x,\mathbf y\in X_F^+$, we define the relation $\sim_+$ as follows:
$$
\mathbf x\sim_+ \mathbf y ~~~   \Longleftrightarrow  ~~~  \mathbf x = \mathbf y \textup{ or  for each positive integer } k,  {\{\mathbf x(k),\mathbf y(k)\}\subseteq} \{0,2,4,6\}. 
$$
\item For all $\mathbf x,\mathbf y\in X_F$, we define the relation $\sim$ as follows:
$$
\mathbf x\sim \mathbf y ~~~   \Longleftrightarrow  ~~~  \mathbf x = \mathbf y \textup{ or  for each integer } k,  {\{\mathbf x(k),\mathbf y(k)\}\subseteq} \{0,2,4,6\}. 
$$
	\end{enumerate}
\end{definition}
\begin{observation}
	Essentially the same proof as the one from \cite[Example 4.14]{BE} shows that the quotient spaces $X_F^+{/}_{\sim_+}$ and $X_F{/}_{\sim}$ are both  Cantor fans. Also, note that $(\sigma_F^+)^{\star}$ is not a homeomorphism on $X_F^+{/}_{\sim_+}$ while $\sigma_F^{\star}$ is a homeomorphism on $X_F{/}_{\sim}$.
	\end{observation}

	\begin{theorem}\label{cantorwsit}
The following hold for the sets of periodic points in $(X_F^+/_{\sim_+},(\sigma_F^+)^{\star})$ and $(X_F/_{\sim},\sigma_F^{\star})$.
\begin{enumerate}
	\item\label{lubi} The set $\mathcal P((\sigma_F^+)^{\star})$ of periodic points in the quotient $(X_F^+/_{\sim_+},(\sigma_F^+)^{\star})$ is dense in $X_F^+/_{\sim_+}$. 
	\item The set $\mathcal P(\sigma_F^{\star})$ of periodic points in the quotient $(X_F/_{\sim},\sigma_F^{\star})$ is dense in $X_F/_{\sim}$. 
\end{enumerate}
	\end{theorem}
\begin{proof}
We prove each of the statements separately.
\begin{enumerate}
	\item We use Theorem \ref{masten} to prove the first part of the theorem. Let $(x,y)\in F$ be any point. We show that there are a positive integer $n$ and a point $\mathbf z\in X_F^n$ such that $\mathbf z(1)=y$ and $\mathbf z(n+1)=x$. We consider the following cases for $x$.
	\begin{enumerate}
	\item $x\in [0,1]$. If $y=x+2$, then let $n=1$ and $\mathbf z=(x+2,x)$. If $y=x^2$, then let $n=3$ and $\mathbf z=(x^2,x^2+2,x+2, x)$. 
	\item $x\in [2,3]$. If $y=x-2$, then let $n=1$ and $\mathbf z=(x-2,x)$. If $y=x+2$, then let $n=1$ and $\mathbf z=(x+2,x)$.  If $y=(x-2)^{\frac{1}{2}}+2$, then let $n=3$ and $\mathbf z=((x-2)^{\frac{1}{2}}+2,(x-2)^{\frac{1}{2}},x-2,x)$. 
	\item $x\in [4,5]$. If $y=x-2$, then let $n=1$ and $\mathbf z=(x-2,x)$. If $y=x+2$, then let $n=1$ and $\mathbf z=(x+2,x)$.  If $y=(x-4)^{3}+4$, then let $n=3$ and $\mathbf z=((x-4)^{3}+4,(x-4)^{3}+6,x+2,x)$. 
	\item $x\in [6,7]$. If $y=x-2$, then let $n=1$ and $\mathbf z=(x-2,x)$. If $y=(x-6)^{\frac{1}{3}}+6$, then let $n=3$ and $\mathbf z=((x-6)^{\frac{1}{3}}+6,(x-6)^{\frac{1}{3}}+4,x-2,x)$.
	\end{enumerate}	
	\item It follows from \ref{lubi} and from Theorem \ref{ingram1} that the set $\mathcal P(\sigma^{-1})$ of periodic points in $(\varprojlim(X_F^+,\sigma_F^+),\sigma^{-1})$ is dense in $\varprojlim(X_F^+,\sigma_F^+)$. By Theorem \ref{povezava}, the set $\mathcal P(\sigma_F^{\star})$ of periodic points in the quotient $(X_F/_{\sim},\sigma_F^{\star})$ is dense in $X_F/_{\sim}$. 
\end{enumerate}
\end{proof}
We use Theorems \ref{A} and \ref{B} to prove that the dynamical systems $(X_F^+/_{\sim_+},(\sigma_F^+)^{\star})$ and $(X_F/_{\sim},\sigma_F^{\star})$ are  transitive.
\begin{definition}
Let $Y$ be a compact metric space,  let $G$ be a closed relation on $Y$ and let $x\in Y$. Then we define 
$$
\mathcal U^{\oplus}_G(x)=\{y\in Y \ | \ \textup{there are } n\in \mathbb N \textup{ and } \mathbf x\in Y_G^{n} \textup{ such that } \mathbf x(1)=x, \mathbf x(n)=y \}
$$
and we call it \emph{the forward impression of $x$ by $G$}.
\end{definition}
\begin{theorem}\label{A}
Let $Y$ be a compact metric space, let $G$ be a closed relation on $Y$,  let  $\{f_{\alpha} \ | \ \alpha \in A\}$ and $\{g_{\beta} \ | \ \beta \in B\}$ be  non-empty collections of continuous functions from $Y$ to $Y$ such that 
$$
G=\bigcup_{\alpha\in A}\Gamma(f_{\alpha}) ~~~ \textup{ and } ~~~ G^{-1}=\bigcup_{\beta\in B}\Gamma(g_{\beta}).
$$
  If there is a dense set $D$ in $Y$ such that for each $s\in D$,  $\Cl(\mathcal U^{\oplus}_G(s))=Y$, then $(Y_G^+,\sigma_G^+)$ is transitive. 
\end{theorem}
\begin{proof}
	See \cite[Theorem 4.8]{BE}.
\end{proof}
\begin{definition}
	For non-empty compact metric spaces $Y$ and $Z$, we use 
$p_1:Y\times Z\rightarrow Y$ and $p_2:Y\times Z\rightarrow Z$
to denote \emph{the standard projections} defined by
$p_1(s,t)=s$ and $p_2(s,t)=t$ for all $(s,t)\in Y\times Z$.
\end{definition}
\begin{theorem}\label{B}
Let $Y$ be a compact metric space and let $G$ be a closed relation on $Y$ such that $p_1(G)=p_2(G)=Y$. The following statements are equivalent. 
\begin{enumerate}
\item The dynamical system  $(Y_G^+,\sigma_G^+)$ is transitive.
\item The dynamical system $(Y_G,\sigma_G)$ is transitive. 
\end{enumerate}
\end{theorem}
\begin{proof}
	The theorem follows from \cite[Theorem 4.5]{BE}.
\end{proof}
We continue with proving results about dynamical systems $(X_F^+/_{\sim_+},(\sigma_F^+)^{\star})$ and $(X_F/_{\sim},\sigma_F^{\star})$ that are defined in Definition \ref{nall2}.
\begin{theorem}\label{transitiveness}
	The dynamical systems $(X_F^+/_{\sim_+},(\sigma_F^+)^{\star})$ and $(X_F/_{\sim},\sigma_F^{\star})$ are both transitive.
	\end{theorem}
\begin{proof}
	To prove that $(X_F^+/_{\sim_+},(\sigma_F^+)^{\star})$ is transitive, we prove that $(X_F^+,\sigma_F^+)$ is transitive.  Note that both $F$ and $F^{-1}$ are unions of three graphs of homeomorphisms. So, all the initial conditions from Theorem \ref{A} are satisfied. To see that $(X_F^+,\sigma_F^+)$ is transitive, we prove that there is a dense set $D$ in $X$ such that for each $s\in D$,  $\Cl(\mathcal U^{\oplus}_H(s))=X$. Let $D=(0,1)\cup(2,3)\cup(4,5)\cup(6,7)$. Then $D$ is dense in $X$. Let $s\in D$ be any point and let $\ell\in\{0,1,2,3\}$ be such that $s\in (2\ell,2\ell+1)$. Note that 
	$$
	s,s-2,s-4, s-6,\ldots, s-2\ell\in \mathcal U^{\oplus}_H(s)
	$$
	and let $t=s-2\ell$. Then $t\in (0,1)$. It follows from the definition of $F$ that for all non-negative integers $m$, $n$ and $k\in\{0,1,2,3\}$, 
	$$
	t^{\frac{2^m}{3^n}}+k\cdot 2\in \mathcal U^{\oplus}_F(t).
	$$
	It follows from Theorem \cite[Lemma 4.9]{BE} that $\big\{t^{\frac{2^m}{3^n}}+k\cdot 2 \ | \ m,n\in \mathbb Z, k\in \{0,1,2,3\}\big\}$ is dense in $X$. Since 
	$$
	\big\{t^{\frac{2^m}{3^n}}+k\cdot 2 \ | \ m,n\in \mathbb Z, k\in \{0,1,2,3\}\big\}\subseteq \mathcal U^{\oplus}_F(t)\subseteq \mathcal U^{\oplus}_F(s),
	$$
	it follows that $\mathcal U^{\oplus}_F(s)$ is dense in $X$. Therefore, by Theorem \ref{A}, $(X_F^+,\sigma_F^+)$ is transitive and it follows from Theorem \ref{B} that $(X_F,\sigma_F)$ is transitive since $p_1(F)=p_2(F)=X$. It follows from Theorem \ref{kvocienti} that the dynamical systems $(X_F^+/_{\sim_+},(\sigma_F^+)^{\star})$ and $(X_F/_{\sim},\sigma_F^{\star})$ are both transitive. 
\end{proof}
\begin{theorem}\label{cantormini}
		The dynamical systems $(X_F^+/_{\sim_+},(\sigma_F^+)^{\star})$ and $(X_F/_{\sim},\sigma_F^{\star})$ both have sensitive dependence on initial conditions. 
	\end{theorem}
\begin{proof}
The dynamical systems $(X_F^+/_{\sim_+},(\sigma_F^+)^{\star})$ and $(X_F/_{\sim},\sigma_F^{\star})$ are both transitive by Theorem \ref{transitiveness}. Also, by Theorem \ref{cantorwsit}, the set $\mathcal P((\sigma_F^+)^{\star})$ of periodic points in the quotient $(X_F^+/_{\sim_+},(\sigma_F^+)^{\star})$ is dense in $X_F^+/_{\sim_+}$,  and the set $\mathcal P(\sigma_F^{\star})$ of periodic points in the quotient $(X_F/_{\sim},\sigma_F^{\star})$ is dense in $X_F/_{\sim}$. It follows from  \cite[Theorem]{banks} that $(X_F^+/_{\sim_+},(\sigma_F^+)^{\star})$ and $(X_F/_{\sim},\sigma_F^{\star})$ both have sensitive dependence on initial conditions.
\end{proof}

\begin{theorem}\label{robi}
The following hold for the Cantor fan $C$. 
\begin{enumerate}
\item There is a continuous mapping $f$ on the Cantor fan $C$ that is not a homeomorphism  such that $(C,f)$ is chaotic in the sense of  Devaney.
	\item There is a homeomorphism $h$ on the Cantor fan $C$ such that $(C,h)$ is chaotic in the sense of Devaney. 
	\end{enumerate}
\end{theorem}
\begin{proof}
We prove each part of the theorem separately. 
\begin{enumerate}
	\item Let $C=X_F^+/_{\sim_+}$ and let $f=(\sigma_F^+)^{\star}$. Note that $f$ is a continuous function which is not a homeomorphism. By Theorem \ref{cantormini}, $(C,f)$ has sensitive dependence on initial conditions; by Theorem \ref{transitiveness}, $(C,f)$ is transitive; and by Theorem \ref{cantorwsit}, the set $\mathcal P(f)$ of periodic points in $(C,f)$ is dense in $C$.  Therefore, $(C,f)$ is chaotic in the sense of Devaney.
	\item Let $C=X_F/_{\sim}$ and let $h=\sigma_F^{\star}$. Note that $h$ is a homeomorphism. By Theorem \ref{cantormini}, $(C,h)$ has sensitive dependence on initial conditions; by Theorem \ref{transitiveness}, $(C,h)$ is transitive; and by Theorem \ref{cantorwsit}, the set $\mathcal P(h)$ of periodic points in  $(C,h)$ is dense in $C$.  Therefore, $(C,h)$ is chaotic in the sense of Devaney.
  \end{enumerate}
 \end{proof}

  \subsection{Robinson's but not Devaney's chaos on the Cantor fan}\label{31} 
          Here, we study functions $f$ on the Cantor fan $C$ such that $(C,f)$ is chaotic in the sense of Robinson but not in the sense of Devaney.  First, we show how Theorem \ref{kvocki} can be used to produce a Cantor fan $C$ and a homeomorphism $h$ on $C$ such that $(C,h)$ has sensitive dependence on initial conditions. 	Our construction of the Cantor fan will use the following definition from \cite{BE}.
	\begin{definition}\label{nall123}
In this subsection, we use $X$ to denote $X=[0,1]\cup [2,3]$. Let 
$$
f_1(x)=\begin{cases}
				x^2\text{;} & x\in [0,1]\\
				(x-2)^{\frac{1}{3}}+2\text{;} & x\in [2,3]
			\end{cases}     ~~~  \textup{ and }  ~~~  f_2(x)=\begin{cases}
				x+2\text{;} & x\in [0,1]\\
				x-2\text{;} & x\in [2,3]
			\end{cases} 
$$
for each $x\in X$. Then we use $F$ to denote the relation $F=\Gamma(f_1)\cup \Gamma(f_2)$; see Figure \ref{fig2}.
\begin{figure}[h!]
	\centering
		\includegraphics[width=13em]{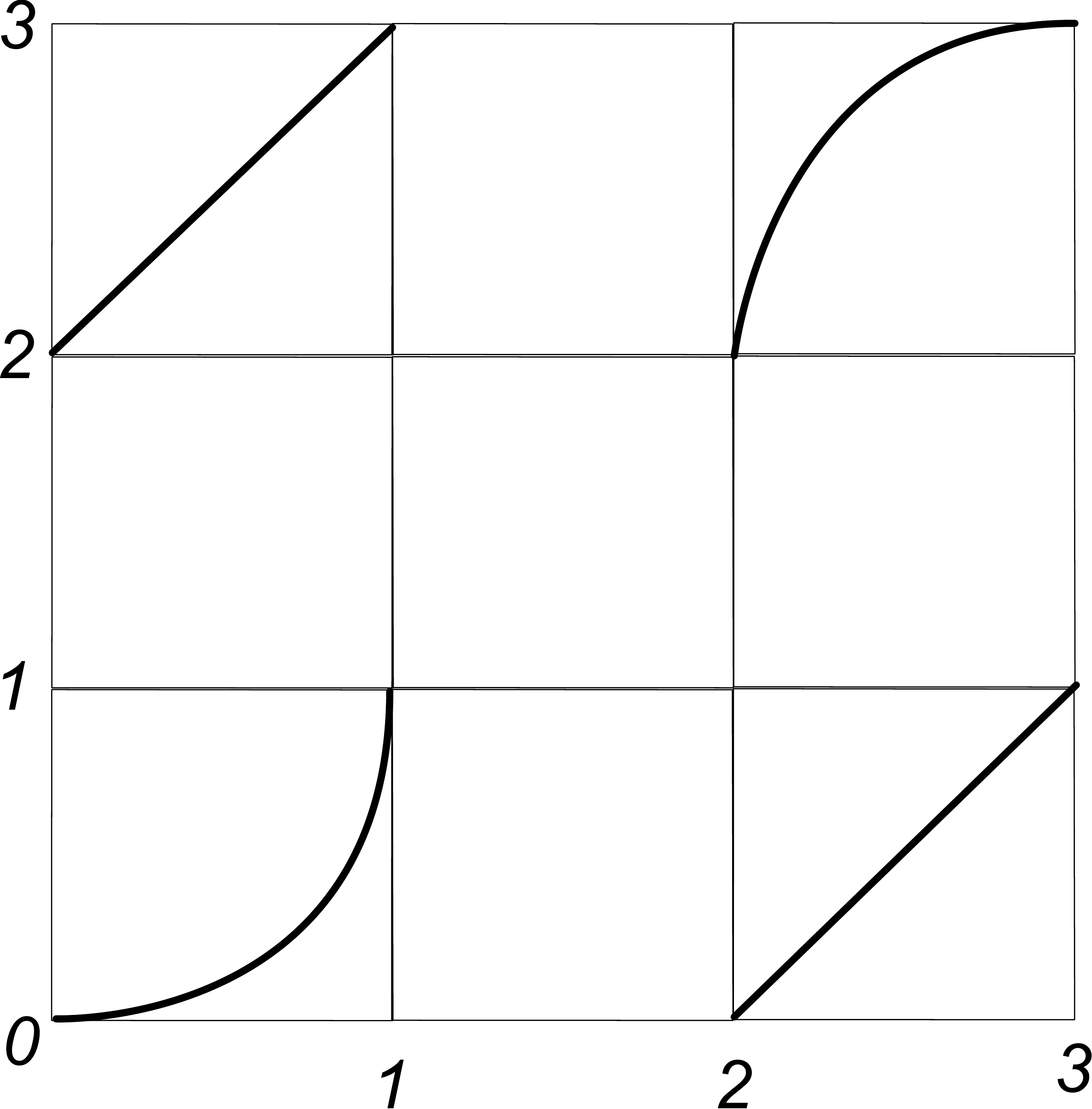}
	\caption{The  relation $F$ from Definition \ref{nall123}}
	\label{fig2}
\end{figure}  
	\end{definition}
	\begin{theorem}\label{skola0}
		Let $A=\big\{\mathbf x\in X_F^+ \ | \ \textup{ for each positive integer } k, ~\mathbf x(k)\in \{0,2\}\big\}$.
		 Then
		 \begin{enumerate}
		 	\item $\sigma_F^+(A)\subseteq A$ and $\sigma_F^+(X_F^+\setminus A)\subseteq X_F^+\setminus A$, and
		 	\item $(X_F^+,\sigma_F^+)$ has sensitive dependence on initial conditions with respect to $A$.
		 \end{enumerate}  
	\end{theorem}
	\begin{proof}
		First, note that $\sigma_F^+(A)\subseteq A$ and $\sigma_F^+(X_F^+\setminus A)\subseteq X_F^+\setminus A$. Next, let $f=\sigma_F^+$ and let $\varepsilon =\frac{1}{4}$. We show that for each basic open set $U$ of the product topology on $\prod_{k=1}^{\infty}X$ such that $U\cap X_F^+\neq \emptyset$, there are $\mathbf x,\mathbf y\in U\cap X_F^+$ such that for some positive integer $m$,
$$
\min\{d(f^m(\mathbf x),f^m(\mathbf y)),d(f^m(\mathbf x),A)+d(f^m(\mathbf y),A)\}>\varepsilon,
$$
 where $d$ is the product metric on $\prod_{k=1}^{\infty}X$, defined defined by  	
		$$
		d((x_1,x_2,x_3,\ldots),(y_1,y_2,y_3,\ldots))=\max\Big\{\frac{|y_k-x_k|}{2^k} \ \big| \  k \textup{ is a positive integer}\Big\}
		$$ 
		for all $(x_1,x_2,x_3,\ldots),(x_1,x_2,x_3,\ldots)\in \prod_{k=1}^{\infty}X$.
  Let $U$ be a basic set of the product topology on $\prod_{k=1}^{\infty}X$ such that $U\cap X_F^+\neq \emptyset$. Also, let $n$ be a positive integer and for each $i\in \{1,2,3,\ldots, n\}$, let $U_i$ be an open set in $X$ such that 
	$$
	U=U_1\times U_2\times U_3\times \ldots \times U_n\times \prod_{k=n+1}^{\infty}X. 
	$$
Next, let $\mathbf z=(z_1,z_2,z_3,\ldots)\in U\cap X_F^+$ be any point such that $z_n\not \in \{0,1,2,3\}$. 
We consider the following possible cases for the coordinate $z_n$ of the point $\mathbf z$.
\begin{enumerate}
	\item $z_n > 2$. Let $\mathbf x=(x_1,x_2,x_3,\ldots), \mathbf y=(y_1,y_2,y_3,\ldots)\in X_F^+$ be such a point that 
	$$
	(x_1,x_2,x_3,\ldots,x_n)=(y_1,y_2,y_3,\ldots,y_n)=(z_1,z_2,z_3,\ldots,z_n)
	$$
	 and for each positive integer $k$, $x_{n+k}$ and $y_{n+k}$ are defined as follows. 	Let  $y_{n+1}=f_2(y_n)$ and for each positive integer $k$, let
	$$
	x_{n+k}=f_1(x_{n+k-1})  ~~~ \textup{ and } ~~~ y_{n+k+1}=f_1(y_{n+k}).
	$$ 
	Note that 
	$$
	\lim_{k\to \infty}f_1^k(x_n)=3 ~~~ \textup{ and } ~~~ \lim_{k\to \infty}f_1^k(f_2(x_n))=0.
	$$
	Let $k_0$ be a positive integer such that for each positive integer $k$,
	$$
	k\geq k_0 ~~~ \Longrightarrow ~~~  3-f_1^k(x_n)<\frac{1}{10} \textup{ and } f_1^k(f_2(x_n))<\frac{1}{10}.
	$$

	Let $m=n+k_0+1$. Then,
	\begin{align*}
		d(f^m(\mathbf x),f^m(\mathbf y))=&\max\Big\{\frac{|y_k-x_k|}{2^{k-m+1}} \ \big| \  k \in \{m,m+1,m+2,m+3,\ldots\}\Big\}=\\
		&\frac{x_m-y_m}{2} = \frac{x_{n+k_0+1}-y_{n+k_0+1}}{2}=\\
		&\frac{f_1^{k_0+1}(x_n)-f_1^{k_0}(f_2(x_n))}{2}\geq \frac{3-\frac{2}{10}}{2}=\frac{7}{5}>\frac{1}{4}=\varepsilon
	\end{align*}
	and 
	\begin{align*}
		d(f^m(\mathbf x),A)+d(f^m(\mathbf y),A)\geq &d(f^m(\mathbf x),A)=\min\{d(f^m(\mathbf x),\mathbf a) \ | \ \mathbf a\in A\}	=\\
		&\min \Bigg\{\max\Big\{\frac{|\mathbf a(k)-x_{k+m}|}{2^{k}} \ \big| \  k \in \{1,2,3,\ldots\}\Big\} \ \Big| \ \mathbf a\in A\Bigg\}\geq \\
		&\min \Bigg\{\max\Big\{\frac{x_{k+m}-\mathbf a(k)}{2^{k}} \ \big| \  k \in \{1,2,3,\ldots\}\Big\} \ \Big| \ \mathbf a\in A\Bigg\}\geq \\
		& \min \Bigg\{\max\Big\{\frac{\big(3-\frac{1}{10}\big)-2}{2^{k}} \ \big| \  k \in \{1,2,3,\ldots\}\Big\}  \ \Big| \ \mathbf a\in A\Bigg\}= \\
		&\min \Bigg\{\max\Big\{\frac{\frac{9}{10}}{2^{k}} \ \big| \  k \in \{1,2,3,\ldots\}\Big\} \ \Big| \ \mathbf a\in A\Bigg\}= \\
		& \max\Big\{\frac{\frac{9}{10}}{2^{k}} \ \big| \  k \in \{1,2,3,\ldots\}\Big\}=\frac{9}{20}>\frac{1}{4}=\varepsilon.
		\end{align*}
	\item $z_n < 1$. In this case, the proof is analogous to the proof of the previous case. We leave the details to the reader.
	    \end{enumerate}
    This proves that $(X_F^+,\sigma_F^+)$ has sensitive dependence on initial conditions with respect to $A$.
	\end{proof}	
	\begin{corollary}\label{corolar}
	Let $B=\big\{\mathbf x\in X_F \ | \ \textup{ for each integer } k, \mathbf x(k)\in \{0,2\}\big\}$.
		 Then
		 \begin{enumerate}
		 	\item $\sigma_F(B)\subseteq B$ and $\sigma_F(X_F\setminus B)\subseteq X_F\setminus B$, and
		 	\item $(X_F,\sigma_F)$ has sensitive dependence on initial conditions with respect to $B$.
		 \end{enumerate}
	\end{corollary}
	\begin{proof}
	First, note that $\sigma_F(B)\subseteq B$ and $\sigma_F(X_F\setminus B)\subseteq X_F\setminus B$. Next, let 
	$$
	A=\big\{\mathbf x\in X_F^+ \ | \ \textup{ for each positive integer } k, \mathbf x(k)\in \{0,2\}\big\}.
	$$
	 By Theorem \ref{skola0},
		 \begin{enumerate}
		 	\item $\sigma_F^+(A)\subseteq A$ and $\sigma_F^+(X_F^+\setminus A)\subseteq X_F^+\setminus A$, and
		 	\item $(X_F^+,\sigma_F^+)$ has sensitive dependence on initial conditions with respect to $A$.
		 \end{enumerate}
		 Note that $\sigma_F^+$ is surjective. By Theorem \ref{ingram}, $		  (\varprojlim(X_F^+,\sigma_F^+),\sigma^{-1})$
		  has sensitive dependence on initial conditions with respect to $\varprojlim(A,\sigma_F^+|_{A})$, where $\sigma$ is the shift homeomorphism on $\varprojlim(X_F^+,\sigma_F^+)$. By Theorem \ref{povezava}, the inverse limit $\varprojlim(X_F^+,\sigma_F^+)$ is homeomorphic to the two-sided Mahavier product $X_F$ and the inverse of the shift homeomorphism $\sigma_F$ on $X_F$ is topologically conjugate to the  shift homeomorphism $\sigma $ on $\varprojlim(X_F^{+},\sigma_F^+)$. Let $ \varphi:\varprojlim(X_F^{+},\sigma_F^+)\rightarrow X_F$ be the homeomorphism, used to prove Theorem \ref{povezava} in \cite[Theorem 4.1]{BE}. Then $\varphi(\varprojlim(A,\sigma_F^+|_{A}))=B$.
 		   	 Therefore, $(X_F,\sigma_F)$ has sensitive dependence on initial conditions with respect to $B$. 
 		   	  \end{proof}
	\begin{definition}
	We define two equivalence relations. 
	\begin{enumerate}
		\item For all $\mathbf x,\mathbf y\in X_F^+$, we define the relation $\sim_+$ as follows:
$$
\mathbf x\sim_+ \mathbf y ~~~   \Longleftrightarrow  ~~~  \mathbf x = \mathbf y \textup{ or  for each positive integer } k,  {\{\mathbf x(k),\mathbf y(k)\}\subseteq} \{0,2\}. 
$$
\item For all $\mathbf x,\mathbf y\in X_F$, we define the relation $\sim$ as follows:
$$
\mathbf x\sim \mathbf y ~~~   \Longleftrightarrow  ~~~  \mathbf x = \mathbf y \textup{ or  for each integer } k,  {\{\mathbf x(k),\mathbf y(k)\}\subseteq} \{0,2\}. 
$$
	\end{enumerate}
\end{definition}
\begin{observation}
	Note that it follows from \cite[Example 4.14]{BE} that the quotient spaces $X_F^+{/}_{\sim_+}$ and $X_F{/}_{\sim}$ are both  Cantor fans. Also, note that $(\sigma_F^+)^{\star}$ is not a homeomorphism on $X_F^+{/}_{\sim_+}$ while $\sigma_F^{\star}$ is a homeomorphism on $X_F{/}_{\sim}$.
	\end{observation}
	\begin{theorem}\label{cantor}
		The dynamical systems $(X_F^+/_{\sim_+},(\sigma_F^+)^{\star})$ and $(X_F/_{\sim},\sigma_F^{\star})$ both have sensitive dependence on initial conditions. 
	\end{theorem}
\begin{proof}
We prove separately for each of the dynamical system that it has sensitive dependence on initial conditions.
\begin{enumerate}
	\item Let $C=X_F^+/_{\sim_+}$ and let $f=(\sigma_F^+)^{\star}$, i.e., for each $\mathbf x\in X_F$, $f([\mathbf x])=[\sigma_F^+(\mathbf x)]$.
We show that $(C,f)$ has sensitive dependence on initial conditions. 
Let 
$$
A=\big\{\mathbf x\in X_F^+ \ | \ \textup{ for each positive integer } k, \mathbf x(k)\in \{0,2\}\big\}.
$$
		 By Theorem  \ref{skola0}, 
		 \begin{enumerate}
		 	\item $\sigma_F(A)\subseteq A$ and $\sigma_F(X_F^+\setminus A)\subseteq X_F^+\setminus A$ and
		 	\item $(X_F^+,\sigma_F^+)$ has sensitive dependence on initial conditions with respect to $A$.
		 \end{enumerate}
Since $A$ is a closed nowhere dense set in $X_F^+$, it follows from Theorem \ref{kvocki} that $(C,f)$ has sensitive dependence on initial conditions.
\item Let $C=X_F/_{\sim}$ and let $h=\sigma_F^{\star}$, i.e., for each $\mathbf x\in X_F$, $h([\mathbf x])=[\sigma_F(\mathbf x)]$.
We show that $(C,h)$ has sensitive dependence on initial conditions. The rest of the proof is analogous to the proof above - instead of the set $A$, the set  
$$
B=\big\{\mathbf x\in X_F \ | \ \textup{ for each integer } k, \mathbf x(k)\in \{0,2\}\big\}
$$
is used in the proof. We leave the details to the reader.
\end{enumerate}
\end{proof}
\begin{theorem}\label{cantorw}
The following hold for the sets of periodic points in $(X_F^+/_{\sim_+},(\sigma_F^+)^{\star})$ and $(X_F/_{\sim},\sigma_F^{\star})$.
\begin{enumerate}
	\item\label{ejn} The set $\mathcal P((\sigma_F^+)^{\star})$ of periodic points in the quotient $(X_F^+/_{\sim_+},(\sigma_F^+)^{\star})$ is not dense in $X_F^+/_{\sim_+}$. 
	\item The set $\mathcal P(\sigma_F^{\star})$ of periodic points in the quotient $(X_F/_{\sim},\sigma_F^{\star})$ is not dense in $X_F/_{\sim}$. 
\end{enumerate}
	\end{theorem}
\begin{proof}
We prove each of the statements separately.
\begin{enumerate}
	\item Suppose that the set $\mathcal P((\sigma_F^+)^{\star})$ of periodic points in  $(X_F^+/_{\sim_+},(\sigma_F^+)^{\star})$ is dense in $X_F^+/_{\sim_+}$. Therefore, by Theorem \ref{kvocki1}, the set $\mathcal P(\sigma_F^+)$ of periodic points in $(X_F^+,\sigma_F^+)$ is dense in $X_F^+$. Let $U=(0,1)\times (0,1)\times \prod_{k=3}^{\infty}X$. Then $U$ is open in $\prod_{k=1}^{\infty}X$ and $U\cap X_F^+\neq \emptyset$. However, note that $U\cap \mathcal P(\sigma_F^+)=\emptyset$ (since for each $x\in (0,1)$, for each $\mathbf x=(x_1,x_2,x_3,\ldots)\in X_F^+$ such that $x_1=x$, and for each positive integer $n>1$, if $x_n\in (0,1)$, then there are positive integers $k$ and $\ell$ such that $x_n=x^{\frac{2^k}{3^{\ell}}}$, which is not equal to $x$). This is a contradiction with the fact that the set $\mathcal P(\sigma_F^+)$ of periodic points in $(X_F^+,\sigma_F^+)$ is dense in $X_F^+$.
	\item Suppose that the set $\mathcal P(\sigma_F^{\star})$ of periodic points in the quotient $(X_F/_{\sim},\sigma_F^{\star})$ is dense in $X_F/_{\sim}$. Therefore, by Theorem \ref{kvocki1}, the set $\mathcal P(\sigma_F)$ of periodic points in $(X_F,\sigma_F)$ is dense in $X_F$. It follows from  Theorem \ref{povezava} that the set $\mathcal P(\sigma^{-1})$ of periodic points in $(\varprojlim(X_F^+,\sigma_F^{+}),\sigma^{-1})$ is dense in $\varprojlim(X_F^+,\sigma_F^{+})$. By Theorem \ref{ingram1}, the set $\mathcal P(\sigma_F^+)$ of periodic points in $(X_F^+,\sigma_F^+)$ is dense in $X_F^+$, which contradicts with \ref{ejn}. 
\end{enumerate}
\end{proof}
\begin{theorem}\label{robi}
The following hold for the Cantor fan $C$. 
\begin{enumerate}
\item There is a continuous mapping $f$ on the Cantor fan $C$ that is not a homeomorphism such that $(C,f)$ is chaotic in the sense of Robinson but not in the sense of Devaney.
	\item There is a homeomorphism $h$ on the Cantor fan $C$ such that $(C,h)$ is chaotic in the sense of Robinson but not in the sense of Devaney. 
	\end{enumerate}
\end{theorem}
\begin{proof}
We prove each part of the theorem separately. 
\begin{enumerate}
	\item Let $C=X_F^+/_{\sim_+}$ and let $f=(\sigma_F^+)^{\star}$. Note that $f$ is a continuous function which is not a homeomorphism. By Theorem \ref{cantor}, $(C,f)$ has sensitive dependence on initial conditions. By \cite[Example 4.14]{BE}, $(C,f)$ is transitive.  It follows from Theorem \ref{cantorw} that the set $\mathcal P(f)$ of periodic points in the quotient $(C,f)$ is not dense in $C$.  Therefore, $(C,f)$ is chaotic in the sense of Robinson but it is not chaotic in the sense of Devaney.
	\item Let $C=X_F/_{\sim}$ and let $h=\sigma_F^{\star}$. Note that $h$ is a homeomorphism. The rest of the proof is analogous to the proof above. We leave the details to the reader.
	 \end{enumerate}
 \end{proof}

     \subsection{Knudsen's but not Devaney's chaos on the Cantor fan}\label{33} 
Here, we study functions $f$ on the Cantor fan $C$ such that $(C,f)$  is chaotic in the sense of Knudsen but not in the sense of Devaney.

\begin{definition}\label{nall1}
In this subsection, we use $X$ to denote $X=[0,1]\cup [2,3]$. Let 
$$
f_1(x)=\begin{cases}
				x^2\text{;} & x\in [0,1]\\
				(x-2)^{\frac{1}{2}}+2\text{;} & x\in [2,3]
			\end{cases}     ~~~  \textup{ and }  ~~~  f_2(x)=\begin{cases}
				x+2\text{;} & x\in [0,1]\\
				x-2\text{;} & x\in [2,3]
			\end{cases} 
$$
for each $x\in X$. Then we use $F$ to denote the relation $F=\Gamma(f_1)\cup \Gamma(f_2)$; see Figure \ref{fig4}.
\begin{figure}[h!]
	\centering
		\includegraphics[width=12em]{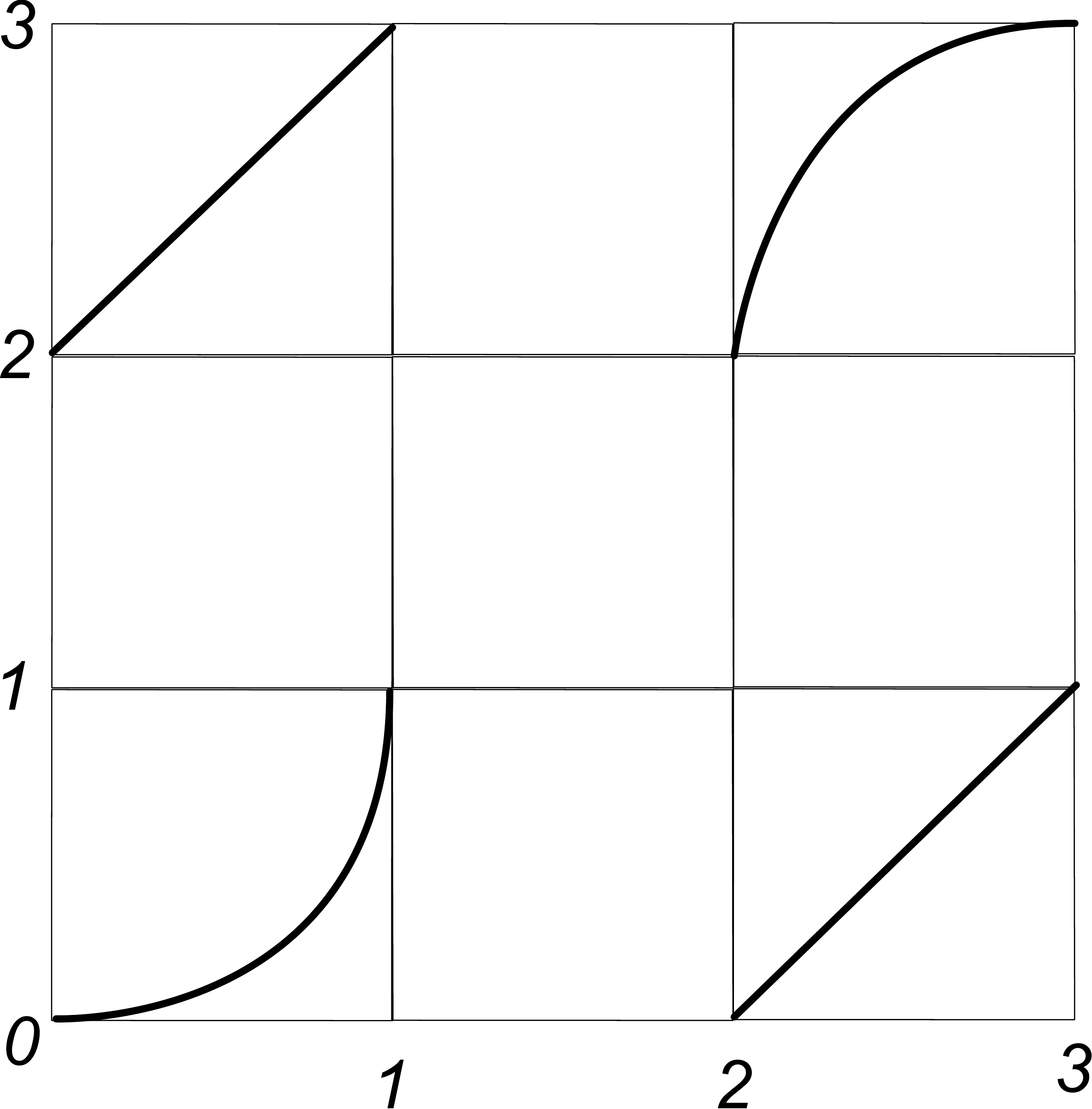}
	\caption{The  relation $F$ from Definition \ref{nall1}}
	\label{fig4}
\end{figure}  
	\end{definition}
	\begin{theorem}\label{skola0s}
		Let $A=\big\{\mathbf x\in X_F^+ \ | \ \textup{ for each positive integer } k, \mathbf x(k)\in \{0,2\}\big\}$.
		 Then
		 \begin{enumerate}
		 	\item $\sigma_F^+(A)\subseteq A$ and $\sigma_F^+(X_F^+\setminus A)\subseteq X_F^+\setminus A$, and
		 	\item $(X_F^+,\sigma_F^+)$ has sensitive dependence on initial conditions with respect to $A$.
		 \end{enumerate}  
	\end{theorem}
	\begin{proof}
		The proof is analogous to the proof of Theorem \ref{skola0}. We leave the details to the reader.
	\end{proof}	
	\begin{corollary}\label{corolars}
	Let $B=\big\{\mathbf x\in X_F \ | \ \textup{ for each integer } k, ~\mathbf x(k)\in \{0,2\}\big\}$.
		 Then
		 \begin{enumerate}
		 	\item $\sigma_F(B)\subseteq B$ and $\sigma_F(X_F\setminus B)\subseteq X_F\setminus B$, and
		 	\item $(X_F,\sigma_F)$ has sensitive dependence on initial conditions with respect to $B$.
		 \end{enumerate}
	\end{corollary}
	\begin{proof}
	The proof is analogous to the proof of Corollary \ref{corolar}. We leave the details to the reader.
 		   	  \end{proof}
	\begin{definition}
	We define two equivalence relations. 
	\begin{enumerate}
		\item For all $\mathbf x,\mathbf y\in X_F^+$, we define the relation $\sim_+$ as follows:
$$
\mathbf x\sim_+ \mathbf y ~~~   \Longleftrightarrow  ~~~  \mathbf x = \mathbf y \textup{ or  for each positive integer } k,  {\{\mathbf x(k),\mathbf y(k)\}\subseteq} \{0,2\}. 
$$
\item For all $\mathbf x,\mathbf y\in X_F$, we define the relation $\sim$ as follows:
$$
\mathbf x\sim \mathbf y ~~~   \Longleftrightarrow  ~~~  \mathbf x = \mathbf y \textup{ or  for each integer } k,  {\{\mathbf x(k),\mathbf y(k)\}\subseteq} \{0,2\}. 
$$
	\end{enumerate}
\end{definition}
\begin{observation}
	Essentially the same proof as the proof of \cite[Example 4.14]{BE} shows that the quotient spaces $X_F^+{/}_{\sim_+}$ and $X_F{/}_{\sim}$ are both  Cantor fans. Again, note that $(\sigma_F^+)^{\star}$ is not a homeomorphism on $X_F^+{/}_{\sim_+}$ while $\sigma_F^{\star}$ is a homeomorphism on $X_F{/}_{\sim}$.
	\end{observation}
	\begin{theorem}\label{cantorss}
		The dynamical systems $(X_F^+/_{\sim_+},(\sigma_F^+)^{\star})$ and $(X_F/_{\sim},\sigma_F^{\star})$ both have sensitive dependence on initial conditions. 
	\end{theorem}
\begin{proof}
The proof is analogous to the proof of Theorem \ref{cantor}. We leave the details to the reader.
\end{proof}
\begin{theorem}\label{cantorwww}
The following hold for the sets of periodic points in $(X_F^+/_{\sim_+},(\sigma_F^+)^{\star})$ and $(X_F/_{\sim},\sigma_F^{\star})$.
\begin{enumerate}
	\item\label{ejn} The set $\mathcal P((\sigma_F^+)^{\star})$ of periodic points in the quotient $(X_F^+/_{\sim_+},(\sigma_F^+)^{\star})$ is dense in $X_F^+/_{\sim_+}$. 
	\item The set $\mathcal P(\sigma_F^{\star})$ of periodic points in the quotient $(X_F/_{\sim},\sigma_F^{\star})$ is dense in $X_F/_{\sim}$. 
\end{enumerate}
	\end{theorem}
\begin{proof}
The proof is analogous to the proof of Theorem \ref{cantorwsit}. We leave the details to the reader.
\end{proof}
\begin{theorem}\label{transtrans}
The following hold for the dynamical systems $(X_F^+/_{\sim_+},(\sigma_F^+)^{\star})$ and $(X_F/_{\sim},\sigma_F^{\star})$.
\begin{enumerate}
	\item\label{ejn} The dynamical system $(X_F^+/_{\sim_+},(\sigma_F^+)^{\star})$ is not transitive. 
	\item The dynamical system $(X_F/_{\sim},\sigma_F^{\star})$ is not transitive.
	\end{enumerate}
	\end{theorem}
\begin{proof}
We prove each of the statements separately.
\begin{enumerate}
	\item To prove that $(X_F^+/_{\sim_+},(\sigma_F^+)^{\star})$ is not transitive, we show first that $(X_F^+,\sigma_F^+)$ is not transitive. Let $x\in X$. We consider the following cases. 
	\begin{enumerate}
	\item $x\in \{0,2\}$. Then $\mathcal U^{\oplus}_F(x)=\{0,2\}$.
	\item $x\in \{1,3\}$. Then $\mathcal U^{\oplus}_F(x)=\{1,3\}$.
	\item $x\not\in \{0,1,2,3\}$. Then 
	$$
	\mathcal U_F^{\oplus}(x)=\{x^{2^k}+2\ell \ | \ k \textup{ is an integer and } \ell\in \{0,1\}\}.
	$$
	\end{enumerate}
	This proves that $\mathcal U_F^{\oplus}(x)$ is not dense in $X$. Let $V$ be a non-empty open set in $X$ such that $V\cap \mathcal U_F^{\oplus}(x)=\emptyset$ and let $U=V\times \prod_{k=2}^{\infty}X$. 	It follows that for each point $\mathbf x=(x_1,x_2,x_3,\ldots)\in X_F^+$ such that $x_1=x$, 
	$$
	\{\mathbf x,\sigma_F^+(\mathbf x), (\sigma_F^+)^2(\mathbf x),(\sigma_F^+)^3(\mathbf x),\ldots\}\cap U=\emptyset.
	$$
	Therefore, for any $\mathbf x\in X_F^+$, the orbit $\{\mathbf x,\sigma_F^+(\mathbf x), (\sigma_F^+)^2(\mathbf x),(\sigma_F^+)^3(\mathbf x),\ldots\}$ of the point $\mathbf x$ is not dense in $X_F^+$. Since $X_F^+$ does not have any isolated points, it follows from Observation \ref{isolatedpoints} that $(X_F^+,\sigma_F^+)$ is not transitive. It follows from Theorem \ref{kvocientek} that the dynamical system $(X_F^+/_{\sim_+},(\sigma_F^+)^{\star})$ is not transitive. 
	\item Since $p_1(F)=p_2(F)=X$ and since $(X_F^+,\sigma_F^+)$ is not transitive, it follows from Theorem \ref{B} that the dynamical system  $(X_F,\sigma_F)$ is not transitive and then, from Theorem \ref{kvocientek} that the dynamical system $(X_F/_{\sim},\sigma_F^{\star})$ is not transitive. 
	\end{enumerate}
\end{proof}

\begin{theorem}\label{robicar}
The following hold for the Cantor fan $C$. 
\begin{enumerate}
\item There is a continuous mapping $f$ on the Cantor fan $C$ that is not a homeomorphism such that $(C,f)$ is chaotic in the sense of Knudsen but not in the sense of Devaney.
	\item There is a homeomorphism $h$ on the Cantor fan $C$ such that $(C,h)$ is chaotic in the sense of Knudsen but not in the sense of Devaney. 
	\end{enumerate}
\end{theorem}
\begin{proof}
We prove each part of the theorem separately. 
\begin{enumerate}
	\item Let $C=X_F^+/_{\sim_+}$ and let $f=(\sigma_F^+)^{\star}$. Note that $f$ is a continuous function which is not a homeomorphism. By Theorem \ref{cantorss}, $(C,f)$ has sensitive dependence on initial conditions. By Theorem \ref{transtrans}, $(C,f)$ is not transitive.  It follows from Theorem \ref{cantorwww} that the set $\mathcal P(f)$ of periodic points in the quotient $(C,f)$ is dense in $C$.  Therefore, $(C,f)$ is chaotic in the sense of Knudsen but it is not chaotic in the sense of Devaney.
	\item Let $C=X_F/_{\sim}$ and let $h=\sigma_F^{\star}$. Note that $h$ is a homeomorphism. The rest of the proof is analogous to the proof above. We leave the details to the reader.
	  \end{enumerate}
 \end{proof}

%%%%%%
%%%%
%%%
%%
%

\section{Chaos on the Lelek fan}\label{s4}
The Lelek fan and the Cantor fan are very different but on the other hand, they share many properties. For example, they both admit transitive homeomorphisms \cite{BE,banic2}, each one of them embeds into the other and they are both universal for the class of smooth fans \cite{universal}.  Therefore, it is only natural to expect that also the Lelek fan might admit similar chaotic properties as does the Cantor fan. So, we finish our paper by stating some open problems about chaos on the Lelek fan. Before that, we give an easy proof (independently from the results about quotients of dynamical systems) that for the Lelek fan $L$,  there is a homeomorphism (as well as a continuous function which is not a homeomorphism) $f:L\rightarrow L$ such hat $(L,f)$ is chaotic in the sense of Robinson but not in the sense of Devaney. 	To do this, we first recall the basic construction from \cite{banic2} of the Lelek fan as the infinite Mahavier product $M_{r,\rho}$ of a closed relation $L_{r,\rho}$ on $[0,1]$.   

\begin{definition}
	Let $X=[0,1]$. For each $(r,\rho)\in (0,\infty)\times (0,\infty)$, we define the sets \emph{$L_r$}, \emph{$L_{\rho}$} and \emph{$L_{r,\rho}$}  as follows:
	\begin{align*}
	&{L_r}=\{(x,y)\in [0,1]\times [0,1] \ | \ y=rx\},
	\\
	&{L_{\rho}}=\{(x,y)\in [0,1]\times [0,1] \ | \ y=\rho x\},
	\\
	&{L_{r,\rho}}=L_r\cup L_{\rho};
	\end{align*}
	see Figure \ref{fig6}.
\begin{figure}[h!]
	\centering
		\includegraphics[width=15em]{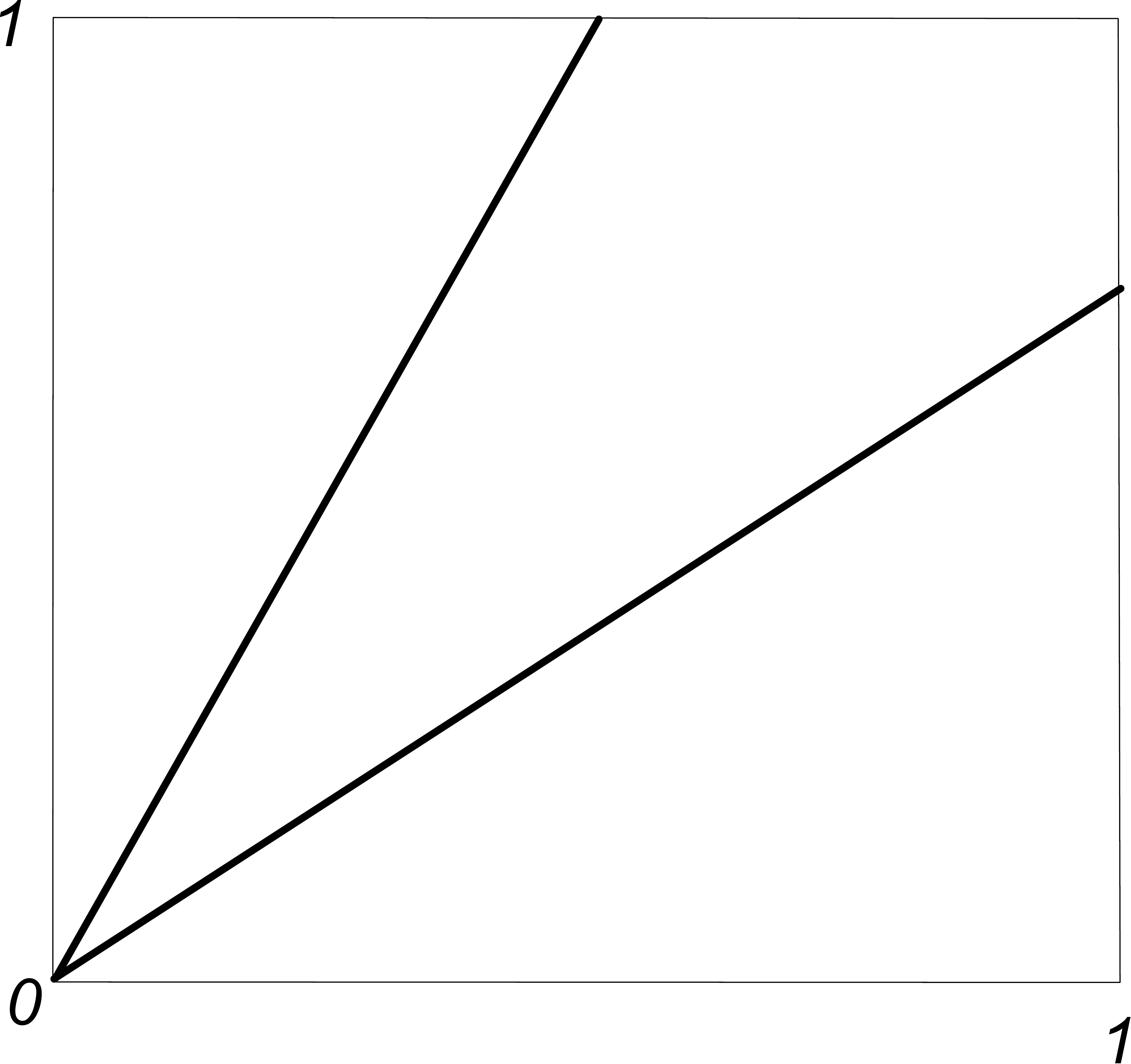}
	\caption{The relation $F$}
	\label{fig6}
\end{figure}  
	
	\noindent We also define the set \emph{$M_{r,\rho}$} by ${M_{r,\rho}}=X_{L_{r,\rho}}^+$ and we use $\sigma_{r,\rho}$ to denote the shift map $\sigma_{L_{r,\rho}}^+$ on $M_{r,\rho}$. 
\end{definition}
\begin{definition}
	Let {$(r,\rho)\in (0,\infty)\times (0,\infty)$}. We say that \emph{$r$ and $\rho$ never connect} or \emph{$(r,\rho)\in \mathcal{NC}$}, if \begin{enumerate}
		\item $r<1$, $\rho>1$ and 
		\item for all integers $k$ and $\ell$,  
		$$
		r^k = \rho^{\ell} \Longleftrightarrow k=\ell=0.
		$$
	\end{enumerate}
	\end{definition} 
	We use the following theorem from \cite{BEK,banic1,banic2}.
	\begin{theorem}\label{Lelek}
Let $(r,\rho)\in \mathcal{NC}$  and let $\mathbf o=(0,0,0,\ldots)$. Then $M_{r,\rho}$ is a Lelek fan with top $\mathbf o$. Furthermore, $(M_{r,\rho},\sigma_{r,\rho})$ is a transitive dynamical system such that $\mathbf o$ is the only periodic point in $(M_{r,\rho},\sigma_{r,\rho})$. 
\end{theorem}
\begin{proof}
It follows from  \cite[Theorem 4.34]{banic1} that $M_{r,\rho}$ is a Lelek fan with top $\mathbf o$. By \cite[Theorem 4.3]{banic2}, $(M_{r,\rho},\sigma_{r,\rho})$ is a transitive dynamical system, and by \cite[Theorem 7.3]{BEK}, $\mathbf o$ is the only periodic point in $(M_{r,\rho},\sigma_{r,\rho})$. 
\end{proof}
In the following theorem, we prove that the dynamical system from Theorem \ref{Lelek} also has sensitive dependence on initial conditions. To prove it, we use the following lemma.
\begin{lemma}\label{goran}
	Let $(r,\rho)\in \mathcal{NC}$. Then for each $x\in (0,1)$,  there is a point 
	$$
	\mathbf x=(x_1,x_2,x_3,\ldots)\in M_{r,\rho}
	$$ 
	such that $x_1=x$ and $\limsup(x_k)=1$. 
	\end{lemma}
\begin{proof}
	The lemma follows from \cite[Theorem 4.26]{banic1}.
\end{proof}
\begin{theorem}\label{SDIC}
Let $(r,\rho)\in \mathcal{NC}$. Then $(M_{r,\rho},\sigma_{r,\rho})$ has sensitive dependence on initial conditions.
\end{theorem}
\begin{proof}
Let $\varepsilon =\frac{1}{4}$. We show that for each basic set $U$ of the product topology on $\prod_{k=1}^{\infty}[0,1]$ such that $U\cap M_{r,\rho}\neq \emptyset$, there are $\mathbf x,\mathbf y\in U\cap M_{r,\rho}$ such that for some positive integer $m$,
$$
d(\sigma_{r,\rho}^m(\mathbf x),\sigma_{r,\rho}^m(\mathbf y))>\varepsilon,
$$
 where $d$ is the product metric on $\prod_{k=1}^{\infty}[0,1]$, defined by  		
		$$
		d((x_1,x_2,x_3,\ldots),(y_1,y_2,y_3,\ldots))=\max\Big\{\frac{|y_k-x_k|}{2^k} \ \big| \  k \textup{ is a positive integer}\Big\}
		$$ 
		for all $(x_1,x_2,x_3,\ldots),(x_1,x_2,x_3,\ldots)\in \prod_{k=1}^{\infty}[0,1]$.
  
  Let $U$ be a basic set of the product topology on $\prod_{k=1}^{\infty}[0,1]$ such that $U\cap M_{r,\rho}\neq \emptyset$. Also, let $n$ be a positive integer and for each $i\in \{1,2,3,\ldots, n\}$, let $U_i$ be an open set in $[0,1]$ such that 
	$$
	U=U_1\times U_2\times U_3\times \ldots \times U_n\times \prod_{k=n+1}^{\infty}[0,1]. 
	$$
Next, let $\mathbf x=(x_1,x_2,x_3,\ldots)\in U\cap M_{r,\rho}$ be any point such that $x_n\not \in \{0,1\}$. 
We observe the following possible cases for the sequence $(x_k)$ of coordinates of $\mathbf x$ in $[0,1]$.
\begin{enumerate}
	\item $\limsup (x_k)>\frac{1}{2}$. Let $(x_{i_k})$ be a subsequence of the sequence $(x_k)$ such that for each positive integer $k$, $x_{i_k}>\frac{1}{2}$, and let 
	$$
	\mathbf y=(y_1,y_2,y_3,\ldots)\in M_{r,\rho}
	$$ 
	be such a point that $(y_1,y_2,y_3,\ldots,y_n)=(x_1,x_2,x_3,\ldots,x_n)$ and for each positive integer $k$, $y_{n+k}=r\cdot y_{n+k-1}$. Since 
	$$
	\lim_{k\to \infty}y_{k}=0 ~~~ \textup{ and for each positive integer } ~~~ k, x_{i_k}>\frac{1}{2},   
 	$$
 	it follows that there is a positive integer $m$ such that $|y_{i_m}-x_{i_m}|>\frac{1}{4}$. Therefore, in this case, there are $\mathbf x,\mathbf y\in U\cap M_{r,\rho}$ such that 
$$
d(\sigma_{r,\rho}^m(\mathbf x),\sigma_{r,\rho}^m(\mathbf y))>\varepsilon,
$$
\item $\limsup(x_k)\leq \frac{1}{2}$. Let $\delta = \frac{1}{10}$ and let $(x_{i_k})$ be a subsequence of the sequence $(x_k)$ such that for each positive integer $k$, $x_{i_k}<\frac{1}{2}+\delta$. Also, let 
	$$
	\mathbf z=(z_1,z_2,z_3,\ldots)\in M_{r,\rho}
	$$ 
	be such a point that $z_1=x_n$ and $\limsup(z_k)=1$. Such a point does exist by Lemma \ref{goran}. Finally, let 	
	$$
	\mathbf y=(y_1,y_2,y_3,\ldots)\in M_{r,\rho}
	$$ 
	be such a point that $(y_1,y_2,y_3,\ldots,y_n)=(x_1,x_2,x_3,\ldots,x_n)$ and for each positive integer $k$, $y_{n+k}=z_{k+1}$. Then $\limsup(y_k)=1$. Since 
	$$
	\limsup(y_k)=1 ~~~ \textup{ and for each positive integer } ~~~ k, x_{i_k}<\frac{1}{2}+\delta,   
 	$$
 	it follows that there is a positive integer $m$ such that $|y_{i_m}-x_{i_m}|>\frac{1}{4}$. Therefore, in this case, there are $\mathbf x,\mathbf y\in U\cap M_{r,\rho}$ such that 
$$
d(\sigma_{r,\rho}^m(\mathbf x),\sigma_{r,\rho}^m(\mathbf y))>\varepsilon.
$$
\end{enumerate}
It follows that there is a positive integer $m$ such that $\diam (\sigma_{r,\rho}^m(U))>\varepsilon$. By Theorem \ref{afna}, $(M_{r,\rho},\sigma_{r,\rho})$ has sensitive dependence on initial conditions.
\end{proof}
Finally, we prove our main result. In its proof, we use the following well-known lemma.
\begin{lemma}\label{iztok}
	Let $(X,f)$ be a dynamical system, where $f:X\rightarrow X$ is a continuous surjection, and let $\sigma$ be the shift homeomorphism on $\varprojlim(X,f)$. If $(X,f)$ has sensitive dependence on initial conditions, then also $(\varprojlim(X,f),\sigma)$ has sensitive dependence on initial conditions. 
\end{lemma}
\begin{proof}
	The proof follows from \cite[Proposition 2.2]{li}.
\end{proof}
\begin{theorem}\label{judy}
	There is a homeomorphism $f:L\rightarrow L$ on the Lelek fan $L$ with top $v$ such that 
	\begin{enumerate}
		\item $(L,f)$ is transitive,
		\item $(L,f)$ has sensitive dependence on initial conditions, and
		\item $\mathcal P(f)=\{v\}$ (so, the set $\mathcal P(f)$ of periodic points in $(L,f)$ is not dense in $L$).
	\end{enumerate}
\end{theorem}
\begin{proof}
Let $(r,\rho)\in \mathcal{NC}$ and let $L=\varprojlim(M_{r,\rho},\sigma_{r,\rho})$. Also, let $v=(\mathbf o, \mathbf o, \mathbf o,\ldots)$, where $\mathbf o=(0,0,0,\ldots)$, and let $f:L\rightarrow L$ be the shift homeomorphism $\sigma$ on $\varprojlim(M_{r,\rho},\sigma_{r,\rho})$. By \cite[Theorem 5.14]{banic2}, $L$ is a Lelek fan. By \cite[Theorem 2.15]{banic2}, $(L,f)$ is transitive. By Theorem \ref{Lelek},  $\mathbf o$ is the only periodic point in $(M_{r,\rho},\sigma_{r,\rho})$. Therefore, $v$ is the only periodic point in $(L,f)$ and it follows that $\mathcal P(f)=\{v\}$ and, therefore, the set $\mathcal P(f)$ of periodic points in $(L,f)$ is not dense in $L$. Finally, it follows from Lemma \ref{iztok} that $(L,f)$ has sensitive dependence on initial conditions. 
\end{proof}
\begin{corollary}
	There is a homeomorphism $h$ on the Lelek fan $L$ such that $(L,h)$ is chaotic in the sense of Robinson but it is not chaotic in the sense of Devaney. Also, there is a continuous mapping $f$ on the Lelek fan $L$, which is not a homeomorphism,  such that $(L,f)$ is chaotic in the sense of Robinson but it is not chaotic in the sense of Devaney.
\end{corollary}
\begin{proof}
	Let $(r,\rho)\in \mathcal{NC}$ and let $L=\varprojlim(M_{r,\rho},\sigma_{r,\rho})$. Also, let $h:L\rightarrow L$ be the shift homeomorphism $\sigma$ on $\varprojlim(M_{r,\rho},\sigma_{r,\rho})$. It follows from the proof of Theorem \ref{judy} that $(L,h)$ is chaotic in the sense of Robinson but it is not chaotic in the sense of Devaney. 
	
	To prove the second part of the corollary, note that it follows from \cite[Theorem 14]{banic1} that $M_{r,\rho}$ is a Lelek fan. Also, note that $\sigma_{r,\rho}:M_{r,\rho}\rightarrow M_{r,\rho}$ is a continuous function which is not a homeomorphism on $M_{r,\rho}$, and that by Theorem \ref{SDIC}, $(M_{r,\rho},\sigma_{r,\rho})$ has sensitive dependence on initial conditions. Also, by \cite[Theorem 4.3]{banic2}, $(M_{r,\rho},\sigma_{r,\rho})$ is transitive.  It follows from Theorem \ref{ingram1} that the set $\mathcal P(\sigma_{r,\rho})$ of periodic points in $(M_{r,\rho},\sigma_{r,\rho})$ is not dense in $M_{r,\rho}$. Therefore, there is a continuous mapping $f$ on the Lelek fan $L$, which is not a homeomorphism,  such that $(L,f)$ is chaotic in the sense of Robinson but it is not chaotic in the sense of Devaney.
\end{proof}
 The following problems are a good place to finish our paper.
\begin{problem}
	Is there a continuous mapping $f$ on the Lelek fan $L$ such that $(L,f)$ is chaotic in the sense of Devaney? 
\end{problem}
\begin{problem}
	Is there a homeomorphism $h$ on the Lelek fan $L$ such that $(L,h)$ is chaotic in the sense of Devaney? 
\end{problem}
\begin{problem}
	Is there a continuous mapping $f$ on the Lelek fan $L$ such that $(L,f)$ is chaotic in the sense of Knudsen but not in the sense of  Devaney? 
\end{problem}
\begin{problem}
	Is there a homeomorphism $h$ on the Lelek fan $L$ such that $(L,h)$ is chaotic in the sense of Knudsen but not in the sense of  Devaney?  
\end{problem}

\section{Acknowledgement}
This work is supported in part by the Slovenian Research Agency (research projects J1-4632, BI-HR/23-24-011, BI-US/22-24-086 and BI-US/22-24-094, and research program P1-0285). 
	
%\section{Declarations}
%The following sections are not relevant to our manuscript. 
%\subsection{Competing interests}
%Not applicable.
%\subsection{Data Availability Statement}
%Not applicable.

\noindent I. Bani\v c\\
              (1) Faculty of Natural Sciences and Mathematics, University of Maribor, Koro\v{s}ka 160, SI-2000 Maribor,
   Slovenia; \\(2) Institute of Mathematics, Physics and Mechanics, Jadranska 19, SI-1000 Ljubljana, 
   Slovenia; \\(3) Andrej Maru\v si\v c Institute, University of Primorska, Muzejski trg 2, SI-6000 Koper,
   Slovenia\\
             {iztok.banic@um.si}           %  \\
     
				\-
				
		\noindent G.  Erceg\\
             Faculty of Science, University of Split, Rudera Bo\v skovi\' ca 33, Split,  Croatia\\
%             {i}     
{{gorerc@pmfst.hr}       }    %  \\

                 	\-
					
  \noindent J.  Kennedy\\
             Department of Mathematics,  Lamar University, 200 Lucas Building, P.O. Box 10047, Beaumont, Texas 77710 USA\\
%             {}     
{{kennedy9905@gmail.com}       }  
\-
				
		\noindent V.  Nall\\
             Department of Mathematics,  University of Richmond, Richmond, Virginia 23173 USA\\
%             {i}     
{{vnall@richmond.edu}       }   

	\-
				

                 %  \\

%``text''
%%%%%%%%%%%%%%%%%%%%%%%%%%%%%%%%%%%%%%%%%%%%%%%%%%%%%%%%%%%%%%%%%%%%%%%%%%%%%%%%%
%%% I N T R O D U C T I O N S

\end{document}